\documentclass[final]{siamltex}
\usepackage{amsmath,color}
\usepackage{amsfonts}
\usepackage{booktabs}
\usepackage[english]{babel}
\usepackage{graphicx}
\usepackage[opta]{optional} 
\usepackage{tikz,caption}
\usepackage{pgfplots}
\pgfplotsset{compat=newest}
\pgfplotsset{plot coordinates/math parser=false}
\usepgfplotslibrary{colormaps} 
\usetikzlibrary{pgfplots.colormaps} 
\usetikzlibrary[pgfplots.colormaps] 
\newlength\figureheight
\newlength\figurewidth 
\usetikzlibrary{external}
\usetikzlibrary{external}
\usepackage[caption=false,font=footnotesize]{subfig}

\usepackage[a4paper]{geometry}

\newtheorem{remark}[theorem]{Remark}
\def\IR{\hbox{\rm I\kern-.2em\hbox{\rm R}}}

\newcommand{\ipmbt}{{\sc ipm-gmres}-${\cal P}_T$ } 
\newcommand{\ipmbd}{{\sc ipm-minres-${\cal P}_D$} }
\newcommand{\ipmbpi}{{\sc ipm-gmres-${\cal P}_\Pi$} }

\newcommand{\ppi}{{\cal P}_\Pi}

\newcommand{\ssnip}{{\sc ssn-gmres-ipf} }

\newcommand{\F}{\mathcal{F}}

\newcommand{\ltwo}{\ensuremath{L^{2}(\Omega)}}

\newcommand{\eps}{\varepsilon}

\newtheorem{teo}{Theorem}[section]

\newcounter{algo}[section]
\renewcommand{\thealgo}{\thesection.\arabic{algo}}
\newcommand{\algo}[3]{\refstepcounter{algo}
\begin{center}\begin{figure}[h]
\framebox[\textwidth]{
\parbox{0.95\textwidth} {\vspace{\topsep}
{\bf Algorithm \thealgo : #2}\label{#1}\\
\vspace*{-\topsep} \mbox{ }\\
{#3} \vspace{\topsep} }}
\end{figure}\end{center}}

\begin{document}
\title{Interior Point Methods and Preconditioning for PDE-Constrained Optimization Problems Involving Sparsity Terms}

\author{John W. Pearson\footnotemark[1] , Margherita Porcelli\footnotemark[2] , and Martin Stoll\footnotemark[3]}
\renewcommand{\thefootnote}{\fnsymbol{footnote}}
\footnotetext[1]{School of Mathematics, The University of Edinburgh, 
        James Clerk Maxwell Building, The King's Buildings, 
        Peter Guthrie Tait Road, Edinburgh, EH9 3FD, United Kingdom
(\texttt{j.pearson@ed.ac.uk})
}
\footnotetext[2]{Universit\`a degli Studi di Firenze, 
Dipartimento di Ingegneria Industriale, 
Viale Morgagni, 40/44, 50134 Firenze, Italy (\texttt{margherita.porcelli@unifi.it})}
\footnotetext[3]{Technische Universit\"{a}t Chemnitz, Faculty
of Mathematics, Professorship Scientific Computing, 09107 Chemnitz, 
Germany
(\texttt{martin.stoll@mathematik.tu-chemnitz.de})
}

\renewcommand{\thefootnote}{\arabic{footnote}}
\maketitle
\begin{abstract}
PDE-constrained optimization problems with control or state constraints are challenging from an analytical as well as numerical perspective. The combination of these constraints with a sparsity-promoting $\rm L^1$ term within the objective function requires sophisticated optimization methods.
We propose the use of an Interior Point scheme applied to a smoothed reformulation
of the discretized problem,
and illustrate that such a scheme exhibits robust performance with respect to parameter changes. 
To increase the potency of this method we introduce fast and efficient preconditioners which enable us to solve problems from a number of PDE applications in low iteration numbers and CPU times, even when the parameters involved are altered dramatically. 
\end{abstract}



\section{Introduction}

In this paper we address the challenge of solving large-scale problems arising from PDE-constrained
optimization \cite{book::hpuu09,book::IK08,book::FT2010}. Such formulations arise in a multitude of applications, ranging
from the control of fluid flows \cite{Hinze2000} to image processing contexts \cite{de2013image}. 
The particular question considered in this paper is how to efficiently handle sparsity-promoting cost
terms within the objective function, as well as additional constraints imposed on the control variable and even the state variable. 
In fact, seeking optimal control functions that are both contained within a range of function
values, and zero on large parts of the domain, has become extremely relevant in practical applications \cite{Sta09}.

In detail, we commence by studying the problem of finding   $(\rm y,\rm u) \in H^1(\Omega) \times L^2(\Omega)$
such that the functional
\begin{align}
\ \F(\rm y,\rm u)&=\frac{1}{2}\|\rm y-\rm y_d\|^ 2_{L^2(\Omega)}+ \frac{\alpha}{2}\|\rm u\|^ 2_{L^2(\Omega)} + 
\beta\|u\|_{L^1(\Omega)}  \label{pb}
\end{align}
is minimized subject to the PDE constraint
\begin{align}
-\Delta \rm y &= \rm u + \rm f ~~\mbox{ in } \Omega, \label{eq:lap} \\        
\rm y  &= \rm g  \hspace{2.3em}\mbox{ on } \Gamma,
\end{align}
where we assume that the equation \eqref{eq:lap} is understood in the weak sense \cite{book::FT2010}. Here, $\Omega\subset\mathbb{R}^2$ or $\mathbb{R}^3$ denotes a spatial domain with boundary $\Gamma$.
Additionally, we allow for box constraints on the control
\begin{equation}\label{box}
\rm u_a \le \rm u \le \rm u_b \quad\mbox{ a.e. in } \Omega,
\end{equation}
and,
for the sake of generality, consider the possibility that there are also box constraints on the state
\begin{equation}\label{boxs}
\rm y_a \le \rm y \le \rm y_b \quad\mbox{ a.e. in } \Omega.
\end{equation}
We follow the convention of recent numerical studies (see \cite{SongYu3,SongYu2,SongYu1,wwL1}, for instance) and investigate the case where the lower (upper) bounds of the box constraints  are non-positive (non-negative). Here, the functions $\rm y_d,f, g,\rm u_a,u_b, y_a,y_b \in \ltwo$ are provided in the problem statement, with $\alpha,\beta >0$ given problem-specific \emph{regularization parameters}. The functions $\rm y,\rm y_d,\rm u$ denote the state, the desired state, and the control, respectively.
The state $\rm y$ and the control $\rm u$ are then linked via a state equation (the PDE). In this work we examine a broad class of state equations, including Poisson's equation (\ref{eq:lap}) as well as the convection--diffusion equation and the heat equation.
Furthermore, we consider the case where the difference between state $\rm y$ and desired state 
$\rm y_d$ is only observed on a certain part of the domain,
i.e. over $\Omega_1\subset\Omega$, with the first quadratic term in (\ref{pb}) 
then having the form $\frac{1}{2}\|\rm y-\rm y_d\|^ 2_{L^2(\Omega_1)}$. We refer to this case
as the ``partial observation'' case.

There are many difficulties associated with the problem (\ref{pb})--(\ref{boxs}), such as selecting a suitable discretization, and choosing an efficient approach for handling the box constraints and the sparsity term. In particular, the state constrained problem itself, not even including the $\rm L^1$-norm term,  leads to a problem formulation where the regularity of the Lagrange multiplier is reduced, see \cite{Cas86} for details. Additionally, the simultaneous treatment of control and state constraints is a complex task.
For this, G\"unther and co-authors in \cite{gunther2012posteriori} propose the use of Moreau--Yosida regularization in order to add the state constraints as a penalty to the objective function. Other approaches are based on a semismooth Newton method, see e.g. \cite{HS10,pss17}.
In fact, the inclusion of control/state constraints leads to a semismooth nonlinear formulation of the first-order optimality 
conditions \cite{BIK99,HIK02,pst15}. Interestingly, the structure of the arising nonlinear system is preserved if the $\rm L^1$-norm 
penalization is added \cite{HS10,pss17,Sta09}. Therefore its solution also generally relies on semismooth Newton approaches, and
an infinite dimensional formulation is commonly utilized to derive the first-order optimality system.  
Stadler in \cite{Sta09} was the first to study PDE-constrained optimization with the  $\rm L^1$ term included,
utilizing a semismooth approach, and  many contributions have been made to the study of these problems in recent years 
(cf.  \cite{HerOW15,HSW11_DS} among others).
Our objective is to tackle the coupled problem of both box constraints combined with the sparsity-promoting term, using the Interior Point method.

The paper \cite{pss17} provides a complete analysis of a globally convergent
semismooth Newton method proposed for the problem (\ref{pb})--(\ref{box}).
Theoretical and practical aspects are investigated 
for both the linear algebra phase and the convergence behavior of the nonlinear method.
The numerical experiments carried out revealed a drawback of the method, as it exhibited 
poor convergence behavior for limiting values of the regularization parameter $\alpha$. 


The aim of this paper is to  propose a new framework for the solution of
(\ref{pb})--(\ref{boxs}) for a wider class of state equations and boundary conditions
and, at the same time, attempt to overcome the numerical limitations of the global semismooth approach.

To pursue this issue we utilize Interior Point methods (IPMs), which 
have shown great applicability for nonlinear programming problems \cite{NocW06,IPMWright}, 
and have also found effective use within the PDE-constrained optimization framework \cite{PGIP17,ulbrich2009primal}.
In particular, IPMs for linear and (convex) quadratic programming problems display 
several features which make them particularly attractive for very large-scale optimization, see e.g. 
the recent survey paper \cite{gondzio12}. Their main advantages are undoubtedly 
their low-degree polynomial worst-case complexity, and 
their ability to deliver optimal solutions in an almost constant number of iterations which
depends very little, if at all, on the problem dimension.
This feature makes IPMs perfect candidates for huge-scale discretized PDE-constrained
optimal control problems.

Recently, in \cite{PGIP17}, an Interior Point approach has been successfully applied to the solution
of problem (\ref{pb})--(\ref{boxs}), with $\beta=0$. In this case the discretization
of the optimization problem  leads to a convex quadratic programming problem,
and IPMs may naturally be applied. Furthermore, the rich structure of the linear systems
arising in this framework allows one to design efficient and robust preconditioners, based on those originally developed for the Poisson control problem without box constraints \cite{PW10}.

In this work we extend the approach proposed in  \cite{PGIP17} to the more difficult 
and general case with $\beta > 0$, and apply it to a broad class of PDE-constrained optimal control problems.
To achieve this goal we utilize two key ingredients that will be described in detail
in Section \ref{sec::ipa}: an appropriate discretization of the
$\rm L^1$-norm that allows us to write the discretized problem in a matrix-vector form, and
a suitable smoothing of the arising vector $\ell_1$-norm that yields a final quadratic programming
form of the discretized problem. The first ingredient is based on the discretization described in  \cite{wwL1},
and recently applied to problem (\ref{pb})--(\ref{box}) in \cite{SongYu3,SongYu2,SongYu1},
where block-coordinate like methods are then introduced.
The second ingredient has been widely used for solving the ubiquitous
$\rm L^1$-norm regularized quadratic problem as, for example, when computing 
sparse solutions in wavelet-based deconvolution problems and compressed sensing  \cite{GPSR}. 
On the other hand, its use is completely new within the PDE-constrained optimization context.
Finally, we propose new preconditioners for the sequence of saddle-point systems
generated by the IPM, based on approximations of the $(1,1)$-block and the Schur complement.
In particular, the case where the $(1,1)$-block is singular is taken into account
when examining the partial observation case.
We may then analyse the spectral properties of the preconditioned $(1,1)$-block and Schur complement, to guide us as to the effectiveness of our overall preconditioning strategies.

We structure the paper as follows. The discretization of the continuous problem is discussed 
in Section \ref{sec::dis}, while an Interior Point scheme is
introduced in Section \ref{sec::ipa} together with the description of the 
linear algebra considerations. Hence, Section \ref{sec::prec} is devoted 
to introducing preconditioning strategies to improve the convergence behavior of the linear iterative solver. 
We highlight a ``matching approach'' that introduces robust approximations to the Schur complement of the linear system.
Additionally, we propose a preconditioning strategy for partial observations in 
Section \ref{subsec::po}, and time-dependent problems in Section \ref{subsec::td}.
Section \ref{exp} illustrates the performance of our scheme for a variety of different parameter regimes, 
discretization levels, and PDE constraints.


\subsection*{Notation}
The $\rm L^1$-norm of a function $\rm u$ is denoted by $\|\rm u\|_{L^1}$,
while the $\ell_1$-norm of a vector $u$ is denoted by $\| u\|_1$. 
Components of a vector $x$ are denoted by $x_j$, or by $x_{a,j}$
for a vector $x_a$. The matrix $I_n$ denotes the $n\times n$ identity matrix,
and $1_n$ is the column vector of ones of dimension $n$.

\section{Problem Discretization and Quadratic Programming Formulation}
\label{sec::dis}
We here apply a discretize-then-optimize approach to (\ref{pb})--(\ref{boxs}), and 
use a finite element discretization that retains a favorable property of the vector $\ell
_1$-norm, specifically that it is separable with respect to the vector components.
This key step allows us to state the discretized problem as a convex quadratic program
that may be tackled using an IPM.

Let $n$ denote the dimension of the discretized space, for both state and control variables. 
Let the matrix $L$ represent a discretization of the Laplacian 
operator (the \textit{stiffness matrix}) when Poisson's equation is considered or, more generally, the discretization of a non-selfadjoint elliptic differential operator, 
and let the matrix $M$ be the finite element Gram matrix, or \textit{mass matrix}.
Finally, we denote by $y,u,y_d,f,u_a,u_b,y_a,y_b$ the discrete counterparts of the functions
$\rm y,u,y_d,f,u_a,u_b,y_a,y_b$, respectively.

The discretization without the additional sparsity term follows a standard Galerkin approach \cite{HS10,RSW09,book::FT2010}.
For the discretization of the $\rm L^1$ term, we here follow \cite{SongYu3,SongYu2,SongYu1,wwL1}
and apply the nodal quadrature rule:
$$\|{\rm u}\|_{{\rm L}^1(\Omega)} \approx \sum^n_{i=1} |u_i| \int_{ \Omega} \phi_i(x)~{\rm d}x,$$ 
where $\{\phi_i\}$ are the finite element basis functions used
and $u_i$ are the components of $u$. It is shown in \cite{wwL1} that first-order convergence may be achieved using this approximation with piecewise linear discretizations of the control. We define a lumped mass matrix $D$ as
$$
D := \text{diag}\left ( \int_{ \Omega} \phi_i(x)~{\rm d}x\right )_{i=1}^{n},
$$
so that the discretized $\rm L^1$-norm can be written in matrix-vector form as $\|D u\|_1$.
As a result, the overall finite element discretization of problem (\ref{pb})--(\ref{boxs}) may be stated as
\begin{equation}
\begin{array}{cl}\label{pb_fe}
\displaystyle\min_{y\in \IR^n,u\in \IR^{n}} & \frac 1 2 (y-y_d)^TM (y-y_d) + 
                          \frac{\alpha}{2}  u^TMu +  \beta \|D u\|_1\\
\mbox{ s.t. }  & L y - Mu = f,
\end{array}
\end{equation}
while additionally being in the presence of control constraints and state constraints: 
\begin{equation}\label{boxvector}
u_a \le  u \le  u_b,\quad\quad y_a \le  y \le  y_b.
\end{equation}
The problems we consider will always have control constraints present, and will sometimes also involve state constraints.

Problem (\ref{pb_fe})--(\ref{boxvector}) is a linearly constrained quadratic problem with bound 
constraints on the state and control variables $(y,u)$, and with an additional nonsmooth weighted 
$\ell_1$-norm term of the variable $u$. 
A possible approach to handle the nonsmoothness in the problem 
consists of using smoothing techniques for the $\ell_1$-norm term, see e.g.  \cite{GPSR,FG-pseudo16,FG-IPM14}.
We here consider a classical strategy proposed in \cite{GPSR} that linearizes
the $\ell_1$-norm by splitting the variable $u$ as follows.
Let $w, v \in \IR^n$ be such that
$$|u_i | = w_i + v_i, \ \ i = 1, \dots, n,
$$
where $w_i = \max(u_i,0)$ and $v_i = \max(-u_i,0)$.  Therefore
$$
\|u\|_1 = 1_n^Tw + 1_n^Tv,  
$$
with $w,v\ge 0$. 
In the weighted case, which we are interested in when approximating the discretized version of $\|\rm u\|_{\rm L^1(\Omega)}$ by $\|Du\|_1$, we obtain 
$$
\|D u\|_1 = 1_n^T Dw + 1_n^T Dv.
$$

By using the relationship 
\begin{equation}\label{split}
 u = w - v,
\end{equation}
one may now rewrite problem (\ref{pb_fe}) in terms of variables $(y,z)$, with 
$$z= \begin{bmatrix}
      w \\
      v
     \end{bmatrix}.
 $$
Note that bounds for $u$
$$u_a \le  u \le  u_b$$
now have to be replaced by the following bounds for $z$:
$$z_a \le  z \le  z_b,$$
with
\begin{equation*}
z_a = \left[\begin{array}{c}
\max\{u_a,0\} \\  -\min\{u_b,0\} \\
\end{array}\right], \qquad  z_b = \left[\begin{array}{c}
\max\{u_b, 0\} \\  -\min\{u_a, 0\} \\
\end{array}\right].
\end{equation*}
We note that these bounds automatically satisfy the constraint $z\ge 0$. Overall, we have the desired quadratic programming formulation:
\begin{equation}
\begin{array}{cl}\label{pb_fe_lin}
\displaystyle\min_{y\in \IR^n,z\in \IR^{2n}} &  Q(y,z):= \frac 1 2 (y-y_d)^TM (y-y_d) + 
                          \frac{\alpha}{2}  z^T \widetilde M z + 
                          \beta\, 1_{2n}^T \bar D z \\
\mbox{ s.t. } & L y - \bar M z = f, \\
               & z_a \le  z \le  z_b, \\
               & y_a \le  y \le  y_b,
\end{array}
\end{equation}
where
$$
 \widetilde M  =  \begin{bmatrix}
                 M & -M \\ 
                 -M & M
                \end{bmatrix}, \quad\quad \bar D  =  \begin{bmatrix}
                 D & D 
                \end{bmatrix} ,\quad\quad
 \bar M   =  \begin{bmatrix}
                 M & -M 
                \end{bmatrix}.
$$
In the next section we derive an Interior Point scheme for the solution of the above problem. 
Clearly once optimal values of variables $z$, and therefore of
$w$ and $v$, are found, the control $u$ of the initial problem is retrieved by
(\ref{split}). We observe that we gain smoothness in the problem at the expense of
increasing  the number of variables by 50\% within the problem statement. 
 Fortunately, this increase will not
have a significant impact in the linear algebra solution phase of our method, as we only require additional sparse matrix-vector multiplications, and the storage of the additional control vectors.

\section{Interior Point Framework and Newton Equations}
\label{sec::ipa}

The three key steps to set up an IPM are the following. First, the
bound constraints are ``eliminated'' by using a logarithmic barrier function.
For problem (\ref{pb_fe_lin}), 
the barrier function takes the form:
\begin{align*}
L_{\mu}(y,z,p) = Q(y,z) + p^T ( L y - \bar M z - f)&{}- \mu \sum \log(y_j - y_{a,j}) - \mu \sum \log(y_{b,j} -y_j)\\
             &{}-\mu \sum \log(z_j - z_{a,j}) - \mu \sum \log(z_{b,j} -z_j),
\end{align*}
where $p\in\IR^n$ is the Lagrange multiplier (or adjoint variable) associated with the state equation,
while $\mu > 0$ is the barrier parameter that controls the relation between
the barrier term and the original objective $Q(y,z)$. As the IPM progresses, $\mu$ is decreased towards zero.

The second step involves applying duality theory, and
deriving the first-order optimality conditions to obtain a nonlinear system
parameterized by $\mu$.
Differentiating $L_\mu$ with respect to $(y,z,p)$ gives the nonlinear system
\begin{eqnarray*}
 M y - M y_d +L^T p - \lambda_{y,a} + \lambda_{y,b} & = & 0, \\
 \alpha \widetilde M z + \beta \bar D ^T 1_{n} - \bar M^T p 
                 - \lambda_{z,a} + \lambda_{z,b} & = & 0, \\
                 L y - \bar M z - f & = & 0,
\end{eqnarray*}
where the $j$th entries of the Lagrange multipliers $\lambda_{y,a},\lambda_{y,b},
\lambda_{z,a},\lambda_{z,b}$ are defined as follows:
$$
(\lambda_{y,a})_j = \frac{\mu}{y_j - y_{a,j}}, \quad\quad
(\lambda_{y,b})_j = \frac{\mu}{y_{b,j} - y_j},  \quad\quad
(\lambda_{z,a})_j = \frac{\mu}{z_j - z_{a,j}},  \quad\quad
(\lambda_{z,b})_j = \frac{\mu}{z_{b,j} - z_j}.
$$
Also, the following bound constraints enforce the 
constraints on $y$ and $z$ via:
$$\lambda_{y,a} \ge 0 , \quad\quad \lambda_{y,b}  \ge 0, \quad\quad  \lambda_{z,a} \ge0, \quad\quad \lambda_{z,b}  \ge  0.$$ 

The third crucial step of the IPM is the application of Newton's method
to the nonlinear system. 
We now derive the Newton equations, following the description in \cite{PGIP17}.
Letting 
$y,z,p, \lambda_{y,a}, \lambda_{y,b}, \lambda_{z,a}, \lambda_{z,b}$
denote the most recent Newton iterates, these quantities 
are updated at each iteration by computing the corresponding Newton steps
$ \Delta y, \Delta z,  \Delta p, \Delta \lambda_{y,a}, \Delta \lambda_{y,b}, \Delta \lambda_{z,a},$ $\Delta \lambda_{z,b}$, 
through the solution of the following Newton system:
\begin{align}
\ \label{7by7} &\begin{bmatrix}
 M             &  0               &   L^T  &  - I_n   &  I_n  &   0  &    0  \\
 0             &  \alpha \widetilde   M   &   -\bar M^T   &  0     &  0   &   -I_{2n} &    I_{2n} \\
 L             &  -\bar M              &    0   &  0     &  0   &   0  &    0 \\
 \Lambda_{y,a} & 0                &    0   & Y - Y_a&  0   &   0  &    0 \\
-\Lambda_{y,b} & 0                &    0   & 0      &Y_b-Y &   0  &    0 \\
0              & \Lambda_{z,a}    &    0   & 0      &0     & Z - Z_a & 0 \\
0              &-\Lambda_{z,b}    &    0   & 0      &0     & 0       & Z_b - Z
\end{bmatrix}
\begin{bmatrix}
 \Delta y\\
 \Delta z \\
 \Delta p \\
 \Delta  \lambda_{y,a} \\
 \Delta  \lambda_{y,b} \\
 \Delta  \lambda_{z,a} \\
 \Delta  \lambda_{z,b} 
\end{bmatrix} \\
\ \nonumber &\hspace{17.5em}=-\begin{bmatrix}
  M y - M y_d +L^T p - \lambda_{y,a} + \lambda_{y,b} \\
 \alpha \widetilde M z + \beta \bar D ^T 1_{n} - \bar M^T p 
                 - \lambda_{z,a} + \lambda_{z,b}  \\
                 L y - \bar M z - f  \\
                 (y-y_a).*\lambda_{y,a} - \mu 1_n \\
                    (y_b-y).*\lambda_{y,b} - \mu 1_n \\
                       (z-z_a).*\lambda_{z,a} - \mu 1_{2n} \\
                          (z_b - z).*\lambda_{z,a} - \mu 1_{2n} 
\end{bmatrix},
\end{align}
where $Y, Z, \Lambda_{y,a}, \Lambda_{y,b}, \Lambda_{z,a}, \Lambda_{z,b}$ are diagonal matrices, 
with the most recent iterates  $y,z,p, \lambda_{y,a},$ $\lambda_{y,b}, \lambda_{z,a}, \lambda_{z,b}$
appearing on their diagonal entries. Similarly, the matrices $Y_a , Y_b , Z_a, Z_b$ are diagonal matrices
corresponding to the bounds $y_a, y_b, z_a, z_b$. 
Here we utilize the {\scshape matlab} notation `$.*$' to denote the componentwise product.
We observe that the contribution of the $\ell_1$-norm term only arises in the right-hand side, that is to say
$\beta$ does not appear within the matrix we need to solve for.

Eliminating $\Delta \lambda_{y,a}, \Delta \lambda_{y,b}, \Delta \lambda_{z,a}, \Delta \lambda_{z,b}$ from \eqref{7by7},
we obtain the following reduced linear system:
\begin{align}\label{NewtonSystem}
&\begin{bmatrix}
  M    + \Theta_y         &  0               &   L^T    \\
 0             &  \alpha \widetilde M  + \Theta_z &   -\bar M^T     \\
 L             &  -\bar M             &    0   \\               
\end{bmatrix}
\begin{bmatrix}
 \Delta y\\
 \Delta z \\
 \Delta p \\
 \end{bmatrix} \\
\nonumber &\hspace{5em}=-\begin{bmatrix}
 M y - M y_d +L^T p  -\mu (Y-Y_a)^{-1}1_n + \mu (Y_b-Y)^{-1}1_n \\
 \alpha \widetilde M z + \beta \bar D ^T 1_{n} - \bar M^T p  
                 -\mu (Z-Z_a)^{-1}1_{2n} + \mu (Z_b-Z)^{-1}1_{2n}  \\
                 L y - \bar M z - f  \\
\end{bmatrix},
\end{align}
with
$$\Theta_y = (Y - Y_a )^{-1} \Lambda_{y,a} + (Y_b - Y )^{-1} \Lambda_{y,b},
\quad\quad\Theta_z = (Z - Z_a )^{-1} \Lambda_{z,a} + (Z_b - Z )^{-1} \Lambda_{z,b}
$$
both diagonal and positive definite matrices, which are typically very ill-conditioned. 
Once the above system is solved, one can compute the steps for the
Lagrange multipliers:
\begin{eqnarray}
 \Delta \lambda_{y,a} & = &  - (Y-Y_a)^{-1} \Lambda_{y,a} \Delta y - \Lambda_{y,a} + \mu (Y-Y_a)^{-1}1_n, \label{zupdate1}\\
 \Delta \lambda_{y,b} & = &   (Y_b-Y)^{-1} \Lambda_{y,b} \Delta y - \Lambda_{y,b} + \mu (Y_b-Y)^{-1}1_n, \label{zupdate2}\\
 \Delta \lambda_{z,a} & = &  - (Z-Z_a)^{-1} \Lambda_{z,a} \Delta z - \Lambda_{z,a} + \mu (Z-Z_a)^{-1}1_{2n}, \label{zupdate3}\\
 \Delta \lambda_{z,b} & = &   (Z_b-Z)^{-1} \Lambda_{z,b} \Delta z - \Lambda_{z,b} + \mu (Z_b-Z)^{-1}1_{2n}. \label{zupdate4}
 \end{eqnarray}
After updating the iterates, and ensuring that they remain feasible, the barrier $\mu$ is reduced and 
a new Newton step is performed.

For the sake of completeness, the structure of the overall Interior Point algorithm is reported in the Appendix, 
and follows the standard infeasible Interior Point path-following scheme outlined in \cite{gondzio12}.
We report on the formulas for the primal and dual feasibilities, given by 
\begin{equation}\label{prdu}
 \xi_p^k =  L y^k - \bar{M} z^k - f, \quad \quad 
 \xi_d^k = \begin{bmatrix}
 M y^k - M y_d + L^T p^k - \lambda^k_{y,a} + \lambda^k_{y,b}   \\
 \alpha \widetilde M z^k + \beta \bar D ^T 1_{n} - \bar {M}^T p^k 
                  - \lambda^k_{z,a} + \lambda^k_{z,b}   
                    \end{bmatrix},
\end{equation}
respectively, and the complementarity gap
\begin{equation}\label{gap}
 \xi_c^k  = \begin{bmatrix}
   (y^k-y_a).* \lambda^k_{y,a} - \mu^k 1_n \\
                    (y_b-y^k).* \lambda^k_{y,b} - \mu^k 1_n \\
                       (z^k-z_a).* \lambda^k_{z,a} - \mu^k 1_{2n} \\
                          (z_b - z^k).* \lambda^k_{z,a} - \mu^k 1_{2n} 
                   \end{bmatrix},
 \end{equation}
for problem (\ref{pb_fe_lin}). Here $k$ denotes the iteration counter for the Interior Point method, with $y^k,z^k,p^k,\lambda^k_{y,a},\lambda^k_{y,b},\lambda^k_{u,a},\lambda^k_{u,b},\mu^k$ the values of $y,z,p,\lambda_{y,a},\lambda_{y,b},\lambda_{u,a},\lambda_{u,b},\mu$ at the $k$th iteration.

The measure of the change in the norm of  $\xi_p^k,  \xi_d^k,  \xi_c^k$ allows us to monitor the convergence of the entire process.
Computationally, the main bottleneck of the algorithm is the linear algebra phase,
that is the efficient solution of the Newton system (\ref{NewtonSystem}).
This is the focus of the forthcoming section.


\section{Preconditioning}
\label{sec::prec}

Having arrived at the Newton system \eqref{NewtonSystem}, the main task at this stage is to construct fast and effective methods for the solution of such systems. In this work, we elect to apply iterative (Krylov subspace) solvers, both the {\scshape minres} method \cite{minres} for symmetric matrix systems, and the {\scshape gmres} algorithm \cite{gmres} which may also be applied to non-symmetric matrices. We wish to accelerate these methods using carefully chosen preconditioners.

To develop these preconditioners, we observe that \eqref{NewtonSystem} is a \emph{saddle-point system} (see \cite{BenGolLie05} for a review of such systems), of the form
\begin{equation*}
\ \mathcal{A}=\left[\begin{array}{cc}
A & B^T \\
B & C \\
\end{array}\right],
\end{equation*}
with
\begin{equation*}
\ A=\left[\begin{array}{cc}
M+\Theta_y & 0 \\
0 & \alpha\widetilde{M}+\Theta_z \\
\end{array}\right],\quad\quad{}B=\left[\begin{array}{cc}
L & -\bar{M} \\
\end{array}\right],\quad\quad{}C=\left[\begin{array}{c}
0 \\
\end{array}\right].
\end{equation*}
Provided $A$ is nonsingular, it is well known that two \emph{ideal preconditioners} for the saddle-point matrix $\mathcal{A}$ are given by
\begin{equation*}
\ \mathcal{P}_1=\left[\begin{array}{cc}
A & 0\\
0 & S \\
\end{array}\right],\quad\quad\mathcal{P}_2=\left[\begin{array}{cc}
A & 0\\
B & -S \\
\end{array}\right],
\end{equation*}
where the (negative) \emph{Schur complement} $S:=-C+BA^{-1}B^T$. In particular, provided the preconditioned system is nonsingular, it can be shown that \cite{Ipsen,Ku95,preconMGW}
\begin{equation*}
\ \lambda(\mathcal{P}_1^{-1}\mathcal{A})\in\left\{1,\frac{1}{2}(1\pm\sqrt{5})\right\},\quad\quad\lambda(\mathcal{P}_2^{-1}\mathcal{A})\in\{1\},
\end{equation*}
and hence that a suitable Krylov method preconditioned by $\mathcal{P}_1$ or $\mathcal{P}_2$ will converge in $3$ or $2$ iterations, respectively.

Of course, we would not wish to work with the preconditioners $\mathcal{P}_1$ or $\mathcal{P}_2$ in practice, as they would be prohibitively expensive to invert. We therefore wish to develop analogous preconditioners of the form
\begin{equation*}
\ \mathcal{P}_D=\left[\begin{array}{cc}
\widehat{A} & 0\\
0 & \widehat{S} \\
\end{array}\right],\quad\quad\mathcal{P}_T=\left[\begin{array}{cc}
\widehat{A} & 0\\
B & -\widehat{S} \\
\end{array}\right],
\end{equation*}
where $\widehat{A}$ and $\widehat{S}$ are suitable and computationally cheap approximations of the $(1,1)$-block $A$ and the Schur complement $S$. Provided $\widehat{A}$ and $\widehat{S}$ are symmetric positive definite, the preconditioner $\mathcal{P}_D$ may be applied within the {\scshape minres} algorithm, and $\mathcal{P}_T$ is applied within a non-symmetric solver such as {\scshape gmres}.

Our focus is therefore to develop such approximations for the corresponding matrices for the Newton system \eqref{NewtonSystem}:
\begin{equation*}
\ A=\left[\begin{array}{cc}
M+\Theta_y & 0 \\
0 & \alpha\widetilde{M}+\Theta_z \\
\end{array}\right],\quad\quad{}S=\left[\begin{array}{cc}
L & -\bar{M} \\
\end{array}\right]\left[\begin{array}{cc}
M+\Theta_y & 0 \\
0 & \alpha\widetilde{M}+\Theta_z \\
\end{array}\right]^{-1}\left[\begin{array}{c}
L^T \\
-\bar{M}^T \\
\end{array}\right].
\end{equation*}

\subsection{Approximation of \boldmath{$(1,1)$}-block}

An effective approximation of the $(1,1)$-block $A$ will require cheap and accurate approximations of the matrices $M+\Theta_y$ and $\alpha\widetilde{M}+\Theta_z$.

When considering the matrix $M+\Theta_y$, our first observation is that the mass matrix $M$ may be effectively approximated by its diagonal \cite{wathen87} within a preconditioner. This can be exploited and enhanced by applying the \emph{Chebyshev semi-iteration} method \cite{VGI61,VGII61,RW08}, which utilizes the effectiveness of the diagonal approximation and accelerates it. Now, it may be easily shown that
\begin{align*}
\ &\Big[\lambda_{\min}\big((D_M+\Theta_y)^{-1}(M+\Theta_y)\big),\lambda_{\max}\big((D_M+\Theta_y)^{-1}(M+\Theta_y)\big)\Big] \\
\ &\hspace{10em}\subset\Big[\min\left\{\lambda_{\min}(D_M^{-1}M),1\right\},\max\left\{\lambda_{\max}(D_M^{-1}M),1\right\}\Big],
\end{align*}
where $D_M:=\text{diag}(M)$, due to the positivity of the diagonal matrix $\Theta_y$. Here, $\lambda_{\min}(\cdot)$, $\lambda_{\max}(\cdot)$ denote the smallest and largest eigenvalues of a matrix, respectively. In other words, the diagonal of $M+\Theta_y$ also clusters the eigenvalues within a preconditioner. The same argument may therefore be used to apply Chebyshev semi-iteration to $M+\Theta_y$ within a preconditioner, and so we elect to use this approach.

We now turn our attention to the matrix $\alpha\widetilde{M}+\Theta_z$, first decomposing $\Theta_z=\text{blkdiag}(\Theta_w,\Theta_v)$, where $\Theta_w$, $\Theta_v$ denote the components of $\Theta_z$ corresponding to $w$, $v$. Therefore, in this notation,
\begin{equation*}
\ \alpha\widetilde{M}+\Theta_z=\left[\begin{array}{cc}
\alpha{}M+\Theta_w & -\alpha{}M \\
-\alpha{}M & \alpha{}M+\Theta_v \\
\end{array}\right].
\end{equation*}
Note that $\widetilde{M}$ is positive semidefinite but $\alpha\widetilde{M}+\Theta_z$ is positive definite since the diagonal $\Theta_z$ is positive definite (the control and state bounds are enforced as strict inequalities at each Newton step).

A result which we apply is that of \cite[Theorems 2.1(i) and 2.2(i)]{LSinverses02}, which gives us the following statements about the inverse of $2\times2$ block matrices:
\begin{teo}
Consider the inverse of the block matrix
\begin{equation}
\ \label{ABCD} \left[\begin{array}{cc}
A & B_1 \\
B_2 & C \\
\end{array}\right].
\end{equation}
If $A$ is nonsingular and $C-B_{2}A^{-1}B_1$ is invertible, then \eqref{ABCD} is invertible, with
\begin{equation}
\ \label{ABCDinv1} \left[\begin{array}{cc}
A & B_1 \\
B_2 & C \\
\end{array}\right]^{-1}=\left[\begin{array}{cc}
A^{-1}+A^{-1}B_1(C-B_{2}A^{-1}B_1)^{-1}B_{2}A^{-1} & -A^{-1}B_1(C-B_{2}A^{-1}B_1)^{-1} \\
-(C-B_{2}A^{-1}B_1)^{-1}B_{2}A^{-1} & (C-B_{2}A^{-1}B_1)^{-1} \\
\end{array}\right].
\end{equation}
Alternatively, if $B_1$ is nonsingular and $B_2-CB_1^{-1}A$ is invertible, then \eqref{ABCD} is invertible, with
\begin{equation}
\ \label{ABCDinv2} \left[\begin{array}{cc}
A & B_1 \\
B_2 & C \\
\end{array}\right]^{-1}=\left[\begin{array}{cc}
-(B_2-CB_1^{-1}A)^{-1}CB_1^{-1} & (B_2-CB_1^{-1}A)^{-1} \\
B_1^{-1}+B_1^{-1}A(B_2-CB_1^{-1}A)^{-1}CB_1^{-1} & -B_1^{-1}A(B_2-CB_1^{-1}A)^{-1} \\
\end{array}\right].
\end{equation}
\end{teo}

For the purposes of this working, we may therefore consider the matrix $\alpha\widetilde{M}+\Theta_z$ itself 
as a block matrix \eqref{ABCD}, 
with $A=\alpha{}M+\Theta_w$, $B_1=B_2=-\alpha{}M$, $C=\alpha{}M+\Theta_v$. It may easily be verified that $A$, $C-B_{2}A^{-1}B_1$, $B_1$, $B_2-CB_1^{-1}A$ are then invertible matrices, and so the results \eqref{ABCDinv1} and \eqref{ABCDinv2} both hold in this setting.

We now consider approximating $\alpha\widetilde{M}+\Theta_z$ within a preconditioner by replacing all mass matrices with their diagonals, i.e. writing
\begin{equation*}
\ \alpha\widetilde{D}_M+\Theta_z:=\left[\begin{array}{cc}
\alpha{}D_M+\Theta_w & -\alpha{}D_M \\
-\alpha{}D_M & \alpha{}D_M+\Theta_v \\
\end{array}\right].
\end{equation*}
This would give us a practical approximation, by using the expression \eqref{ABCDinv1} to apply $(\alpha\widetilde{D}_M+\Theta_z)^{-1}$, provided it can be demonstrated that $\alpha\widetilde{D}_M+\Theta_z$ well approximates $\alpha\widetilde{M}+\Theta_z$. This is indeed the case, as demonstrated using the result below:
\begin{teo}
	\label{theorem1}
The eigenvalues $\lambda$ of the matrix
\begin{equation}
\ \label{PrecDiag} \left[\begin{array}{cc}
\alpha{}D_M+\Theta_w & -\alpha{}D_M \\
-\alpha{}D_M & \alpha{}D_M+\Theta_v \\
\end{array}\right]^{-1}\left[\begin{array}{cc}
\alpha{}M+\Theta_w & -\alpha{}M \\
-\alpha{}M & \alpha{}M+\Theta_v \\
\end{array}\right]
\end{equation}
are all contained within the interval:
\begin{equation*}
\ \lambda\in\Big[\min\{\lambda_{\min}(D_M^{-1}M),1\},\max\{\lambda_{\max}(D_M^{-1}M),1\}\Big].
\end{equation*}
\end{teo}
\emph{Proof.}~~The eigenvalues of \eqref{PrecDiag} satisfy
\begin{equation*}
\ \left[\begin{array}{cc}
\alpha{}M+\Theta_w & -\alpha{}M \\
-\alpha{}M & \alpha{}M+\Theta_v \\
\end{array}\right]\left[\begin{array}{c}
\mathbf{x}_1 \\
\mathbf{x}_2 \\
\end{array}\right]=\lambda\left[\begin{array}{cc}
\alpha{}D_M+\Theta_w & -\alpha{}D_M \\
-\alpha{}D_M & \alpha{}D_M+\Theta_v \\
\end{array}\right]\left[\begin{array}{c}
\mathbf{x}_1 \\
\mathbf{x}_2 \\
\end{array}\right],
\end{equation*}
with $\mathbf{x}_1$, $\mathbf{x}_2$ not both equal to $\mathbf{0}$, which may be decomposed to write
\begin{align}
\ \label{EigEqn1} (\alpha{}M+\Theta_w)\mathbf{x}_1-\alpha{}M\mathbf{x}_2={}&\lambda(\alpha{}D_M+\Theta_w)\mathbf{x}_1-\lambda\alpha{}D_M\mathbf{x}_2, \\
\ \label{EigEqn2} -\alpha{}M\mathbf{x}_1+(\alpha{}M+\Theta_v)\mathbf{x}_2={}&-\lambda\alpha{}D_M\mathbf{x}_1+\lambda(\alpha{}D_M+\Theta_v)\mathbf{x}_2.
\end{align}
Summing \eqref{EigEqn1} and \eqref{EigEqn2} gives that
\begin{equation*}
\ \Theta_w\mathbf{x}_1+\Theta_v\mathbf{x}_2=\lambda{}\Theta_w\mathbf{x}_1+\lambda{}\Theta_v\mathbf{x}_2=\lambda(\Theta_w\mathbf{x}_1+\Theta_v\mathbf{x}_2),
\end{equation*}
which tells us that either $\lambda=1$ or $\Theta_w\mathbf{x}_1+\Theta_v\mathbf{x}_2=\mathbf{0}$. In the latter case, we substitute $\mathbf{x}_1=-\Theta_w^{-1}\Theta_v\mathbf{x}_2$ into \eqref{EigEqn1} to give that
\begin{align*}
\ -(\alpha{}M+\Theta_w)\Theta_w^{-1}\Theta_v\mathbf{x}_2-\alpha{}M\mathbf{x}_2={}&-\lambda(\alpha{}D_M+\Theta_w)\Theta_w^{-1}\Theta_v\mathbf{x}_2-\lambda\alpha{}D_M\mathbf{x}_2 \\
\ \Rightarrow\quad\quad~~\Big[\alpha{}M(\Theta_w^{-1}\Theta_v+I)+\Theta_v\Big]\mathbf{x}_2={}&\lambda\Big[\alpha{}D_M(\Theta_w^{-1}\Theta_v+I)+\Theta_v\Big]\mathbf{x}_2,
\end{align*}
which in turn tells us that
\begin{equation*}
\ \Big[\alpha{}M(\Theta_w^{-1}\Theta_v+I)^{1/2}+\Theta_v(\Theta_w^{-1}\Theta_v+I)^{-1/2}\Big]\mathbf{x}_3=\lambda\Big[\alpha{}D_M(\Theta_w^{-1}\Theta_v+I)^{1/2}+\Theta_v(\Theta_w^{-1}\Theta_v+I)^{-1/2}\Big]\mathbf{x}_3,
\end{equation*}
where $\mathbf{x}_3=(\Theta_w^{-1}\Theta_v+I)^{1/2}\mathbf{x}_2\neq\mathbf{0}$. Premultiplying both sides of the equation by $(\Theta_w^{-1}\Theta_v+I)^{1/2}$ then gives that
\begin{equation*}
\ \Big[\alpha(\Theta_w^{-1}\Theta_v+I)^{1/2}M(\Theta_w^{-1}\Theta_v+I)^{1/2}+\Theta_v\Big]\mathbf{x}_3=\lambda\Big[\alpha(\Theta_w^{-1}\Theta_v+I)^{1/2}D_M(\Theta_w^{-1}\Theta_v+I)^{1/2}+\Theta_v\Big]\mathbf{x}_3,
\end{equation*}
and therefore that the eigenvalues may be described by the Rayleigh quotient
\begin{equation*}
\ \frac{\mathbf{x}_3^T\Big[\alpha(\Theta_w^{-1}\Theta_v+I)^{1/2}M(\Theta_w^{-1}\Theta_v+I)^{1/2}+\Theta_v\Big]\mathbf{x}_3}{\mathbf{x}_3^T\Big[\alpha(\Theta_w^{-1}\Theta_v+I)^{1/2}D_M(\Theta_w^{-1}\Theta_v+I)^{1/2}+\Theta_v\Big]\mathbf{x}_3}.
\end{equation*}
Now, as $\mathbf{x}_3^T\Theta_v\mathbf{x}_3$ is a positive number, $\lambda$ may be bounded within the range of the following Rayleigh quotient:
\begin{align*}
\ \lambda\in{}&\left[\min\left\{\min_{\mathbf{x}_3}\frac{\mathbf{x}_3^T\Big[\alpha(\Theta_w^{-1}\Theta_v+I)^{1/2}M(\Theta_w^{-1}\Theta_v+I)^{1/2}\Big]\mathbf{x}_3}{\mathbf{x}_3^T\Big[\alpha(\Theta_w^{-1}\Theta_v+I)^{1/2}D_M(\Theta_w^{-1}\Theta_v+I)^{1/2}\Big]\mathbf{x}_3},1\right\},\right. \\
\ &\quad\quad\quad\quad\left.\max\left\{\max_{\mathbf{x}_3}\frac{\mathbf{x}_3^T\Big[\alpha(\Theta_w^{-1}\Theta_v+I)^{1/2}M(\Theta_w^{-1}\Theta_v+I)^{1/2}\Big]\mathbf{x}_3}{\mathbf{x}_3^T\Big[\alpha(\Theta_w^{-1}\Theta_v+I)^{1/2}D_M(\Theta_w^{-1}\Theta_v+I)^{1/2}\Big]\mathbf{x}_3},1\right\}\right] \\
\ ={}&\left[\min\left\{\min_{\mathbf{x}_4}\frac{\mathbf{x}_4^{T}M\mathbf{x}_4}{\mathbf{x}_4^{T}D_M\mathbf{x}_4},1\right\},\max\left\{\max_{\mathbf{x}_4}\frac{\mathbf{x}_4^{T}M\mathbf{x}_4}{\mathbf{x}_4^{T}D_M\mathbf{x}_4},1\right\}\right] \\
\ \subset{}&\Big[\min\{\lambda_{\min}(D_M^{-1}M),1\},\max\{\lambda_{\max}(D_M^{-1}M),1\}\Big],
\end{align*}
where in the above derivation $\mathbf{x}_4=(\Theta_w^{-1}\Theta_v+I)^{1/2}\mathbf{x}_3\neq\mathbf{0}$. This gives the stated result.~~$\Box$

\vspace{1em}

\begin{remark} Theorem \ref{theorem1} is indeed a positive result. We utilize the fact that a mass matrix preconditioned by its diagonal gives tight eigenvalue bounds \cite{wathen87}. 
We have now obtained a cheap approximation of the $(1,1)$-block of our saddle-point system, with eigenvalues of the preconditioned matrix provably contained within a tight interval. 
We wish to emphasize the fact that the interval boundaries and thus the region of interest where the eigenvalues will lie is independent of all system parameters, such as penalization-, regularization-, mesh-, and time-step parameters.
\end{remark}

\subsection{Approximation of Schur Complement}\label{sec:Schur}

The Schur complement of the Newton system \eqref{NewtonSystem} under consideration is given by
\begin{equation*}
\ S=L(M+\Theta_y)^{-1}L^{T}+\left[\begin{array}{cc}
-M & M \\
\end{array}\right]\left[\begin{array}{cc}
\alpha{}M+\Theta_w & -\alpha{}M \\
-\alpha{}M & \alpha{}M+\Theta_v \\
\end{array}\right]^{-1}\left[\begin{array}{c}
-M \\
M \\
\end{array}\right].
\end{equation*}
For the matrix inverse in the above expression, we again consider the matrix $\alpha\widetilde{M}+\Theta_z$ as
a block matrix of the form \eqref{ABCD}, with $A=\alpha{}M+\Theta_w$, $B_1=B_2=B=-\alpha{}M$, $C=\alpha{}M+\Theta_v$. Using \eqref{ABCDinv2} then gives that
\begin{align*}
\ &\left[\begin{array}{cc}
-M & M \\
\end{array}\right]\left[\begin{array}{cc}
A & B \\
B & C \\
\end{array}\right]^{-1}\left[\begin{array}{c}
-M \\
M \\
\end{array}\right] \\
\ ={}&\left[\begin{array}{cc}
-M & M \\
\end{array}\right]\left[\begin{array}{c}
(B-CB^{-1}A)^{-1}CB^{-1}M+(B-CB^{-1}A)^{-1}M \\
-B^{-1}M-B^{-1}A(B-CB^{-1}A)^{-1}CB^{-1}M-B^{-1}A(B-CB^{-1}A)^{-1}M \\
\end{array}\right] \\
\ \ ={}&-M\Big[B^{-1}+(B^{-1}A+I)(B-CB^{-1}A)^{-1}(CB^{-1}+I)\Big]M,
\end{align*}
whereupon substituting in the relevant $A$, $B$, $C$ gives that this expression can be written as follows:
\begin{align*}
\ &\frac{1}{\alpha}M-\left(-\frac{1}{\alpha}A+M\right)\left(-\alpha{}M+\frac{1}{\alpha}DM^{-1}A\right)^{-1}\left(-\frac{1}{\alpha}D+M\right) \\
\ ={}&\frac{1}{\alpha}M+\left(\frac{1}{\alpha}\Theta_w\right)\left(\alpha{}M-\left(\alpha{}M+\Theta_w+\Theta_v+\frac{1}{\alpha}\Theta_v{}M^{-1}\Theta_w\right)\right)^{-1}\left(\frac{1}{\alpha}\Theta_v\right) \\
\ ={}&\frac{1}{\alpha}M-\frac{1}{\alpha^2}\left(\Theta_w^{-1}+\Theta_v^{-1}+\frac{1}{\alpha}M^{-1}\right)^{-1}.
\end{align*}
Therefore, $S$ may be written as
\begin{equation}
\ \label{Schur} S=L(M+\Theta_y)^{-1}L^{T}+\frac{1}{\alpha}M-\frac{1}{\alpha^2}\left(\Theta_w^{-1}+\Theta_v^{-1}+\frac{1}{\alpha}M^{-1}\right)^{-1}.
\end{equation}
It can be shown that $S$ consists of a sum of two symmetric positive semidefinite matrices. The matrix $L(M+\Theta_y)^{-1}L^{T}$ clearly satisfies this property due to the positive definiteness of $M+\Theta_y$, and $\frac{1}{\alpha}M-\frac{1}{\alpha^2}\left(\Theta_w^{-1}+\Theta_v^{-1}+\frac{1}{\alpha}M^{-1}\right)^{-1}$ is in fact positive definite by the following argument:
\begin{align*}
\frac{1}{\alpha}M-\frac{1}{\alpha^2}\left(\frac{1}{\alpha}M^{-1}+\Theta_w^{-1}+\Theta_v^{-1}\right)^{-1}\succ0\quad\Leftrightarrow\quad&\frac{1}{\alpha^2}\left(\frac{1}{\alpha}M^{-1}+\Theta_w^{-1}+\Theta_v^{-1}\right)^{-1}\prec\frac{1}{\alpha}M \\
\ \Leftrightarrow\quad&\alpha^2\left(\frac{1}{\alpha}M^{-1}+\Theta_w^{-1}+\Theta_v^{-1}\right)\succ\alpha{}M^{-1} \\
\ \Leftrightarrow\quad&\ M^{-1}+\alpha\Theta_w^{-1}+\alpha\Theta_v^{-1}\succ{}M^{-1}.
\end{align*}
Based on this observation, we apply a ``\emph{matching strategy}'' previously derived in \cite{PSW11,PW10} for simpler PDE-constrained optimization problems, which relies on a Schur complement being written in this form. In more detail, we approximate the Schur complement $S$ by
\begin{equation}
\ \label{SchurApprox} \widehat{S}=\left(L+\widehat{M}\right)(M+\Theta_y)^{-1}\left(L+\widehat{M}\right)^T,
\end{equation}
where $\widehat{M}$ is chosen such that the `outer' term of $\widehat{S}$ in \eqref{SchurApprox} approximates the second and third terms of $S$ in \eqref{Schur}, that is
\begin{equation*}
\ \widehat{M}(M+\Theta_y)^{-1}\widehat{M}^{T}\approx\frac{1}{\alpha}M-\frac{1}{\alpha^2}\left(\Theta_w^{-1}+\Theta_v^{-1}+\frac{1}{\alpha}M^{-1}\right)^{-1}.
\end{equation*}
This may be achieved if
\begin{equation*}
\ \widehat{M}\approx\left[\frac{1}{\alpha}M-\frac{1}{\alpha^2}\left(\Theta_w^{-1}+\Theta_v^{-1}+\frac{1}{\alpha}M^{-1}\right)^{-1}\right]^{1/2}(M+\Theta_y)^{1/2}.
\end{equation*}
A natural choice, which may be readily worked with on a computer, therefore involves replacing mass matrices with their diagonals, making the square roots of matrices practical to work with, and therefore setting
\begin{equation*}
\ \widehat{M}=\left[\frac{1}{\alpha}D_M-\frac{1}{\alpha^2}\left(\Theta_w^{-1}+\Theta_v^{-1}+\frac{1}{\alpha}D_M^{-1}\right)^{-1}\right]^{1/2}(D_M+\Theta_y)^{1/2}.
\end{equation*}
We therefore have a Schur complement approximation $\widehat{S}$ which may be approximately inverted by applying a multigrid method to the matrix $L+\widehat{M}$ and its transpose, along with a matrix-vector multiplication for $M+\Theta_y$.

Below we present a result concerning the lower bounds of the eigenvalues of the preconditioned Schur complement.
\begin{teo}
In the case of lumped (diagonal) mass matrices, the eigenvalues of the preconditioned Schur complement all satisfy:
\begin{equation*}
\ \lambda(\widehat{S}^{-1}S)\geq\frac{1}{2}.
\end{equation*}
\end{teo}
\emph{Proof.}~~Bounds for the eigenvalues of $\widehat{S}^{-1}S$ are determined by the extrema of the Rayleigh quotient
\begin{equation*}
\ R:=\frac{\mathbf{v}^{T}S\mathbf{v}}{\mathbf{v}^{T}\widehat{S}\mathbf{v}}=\frac{\boldsymbol\chi^T\boldsymbol\chi+\boldsymbol\omega^T\boldsymbol\omega}{(\boldsymbol\chi+\boldsymbol\gamma)^T(\boldsymbol\chi+\boldsymbol\gamma)},
\end{equation*}
where
\begin{align*}
\ \boldsymbol\chi={}&(M+\Theta_y)^{-1/2}L^T\mathbf{v}, \\
\ \boldsymbol\omega={}&\left[\frac{1}{\alpha}M-\frac{1}{\alpha^2}\left(\Theta_w^{-1}+\Theta_v^{-1}+\frac{1}{\alpha}M^{-1}\right)^{-1}\right]^{1/2}\mathbf{v}, \\
\ \boldsymbol\gamma={}&(M+\Theta_y)^{-1/2}(D_M+\Theta_y)^{1/2}\left[\frac{1}{\alpha}D_M-\frac{1}{\alpha^2}\left(\Theta_w^{-1}+\Theta_v^{-1}+\frac{1}{\alpha}D_M^{-1}\right)^{-1}\right]^{1/2}\mathbf{v}.
\end{align*}
Following the argument used in \cite[Lemma 2]{PGIP17}, we may bound $R$ as follows:
\begin{equation}\label{Rbound}
\ R=\frac{\boldsymbol\chi^T\boldsymbol\chi+\displaystyle{\frac{\boldsymbol\omega^T\boldsymbol\omega}{\boldsymbol\gamma^T\boldsymbol\gamma}}\hspace{0.25em}\boldsymbol\gamma^T\boldsymbol\gamma}{(\boldsymbol\chi+\boldsymbol\gamma)^T(\boldsymbol\chi+\boldsymbol\gamma)}\geq\min\left\{\frac{\boldsymbol\omega^T\boldsymbol\omega}{\boldsymbol\gamma^T\boldsymbol\gamma},1\right\}\cdot\frac{\boldsymbol\chi^T\boldsymbol\chi+\boldsymbol\gamma^T\boldsymbol\gamma}{(\boldsymbol\chi+\boldsymbol\gamma)^T(\boldsymbol\chi+\boldsymbol\gamma)}\geq\frac{1}{2}\cdot\min\left\{\frac{\boldsymbol\omega^T\boldsymbol\omega}{\boldsymbol\gamma^T\boldsymbol\gamma},1\right\},
\end{equation}
using the argument
\begin{align*}
\ \frac{1}{2}(\boldsymbol\chi-\boldsymbol\gamma)^T(\boldsymbol\chi-\boldsymbol\gamma)\geq0\quad\Leftrightarrow&\quad\boldsymbol\chi^T\boldsymbol\chi+\boldsymbol\gamma^T\boldsymbol\gamma\geq\frac{1}{2}(\boldsymbol\chi+\boldsymbol\gamma)^T(\boldsymbol\chi+\boldsymbol\gamma) \\
\ \Leftrightarrow&\quad\frac{\boldsymbol\chi^T\boldsymbol\chi+\boldsymbol\gamma^T\boldsymbol\gamma}{(\boldsymbol\chi+\boldsymbol\gamma)^T(\boldsymbol\chi+\boldsymbol\gamma)}\geq\frac{1}{2}.
\end{align*}

We now turn our attention to the product $\frac{\boldsymbol\omega^T\boldsymbol\omega}{\boldsymbol\gamma^T\boldsymbol\gamma}$. Straightforward calculation tells us that
\begin{equation*}
\ \frac{\boldsymbol\omega^T\boldsymbol\omega}{\boldsymbol\gamma^T\boldsymbol\gamma}=\underbrace{\frac{\mathbf{v}^T[M-(\Theta+M^{-1})^{-1}]\mathbf{v}}{\mathbf{v}^T[D_M-(\Theta+D_M^{-1})^{-1}]\mathbf{v}}}_{=:R_{\Theta}}\cdot\frac{\mathbf{w}^T(D_M+\Theta_y)^{-1}\mathbf{w}}{\mathbf{w}^T(M+\Theta_y)^{-1}\mathbf{w}},
\end{equation*}
where $\Theta:=\alpha\Theta_w^{-1}+\alpha\Theta_v^{-1}$ and $\mathbf{w}:=(D_M+\Theta_y)^{1/2}\big[\frac{1}{\alpha}D_M-\frac{1}{\alpha^2}\left(\Theta_w^{-1}+\Theta_v^{-1}+\frac{1}{\alpha}D_M^{-1}\right)^{-1}\big]^{1/2}\mathbf{v}$. It may be observed that
\begin{equation*}
\ \frac{\mathbf{w}^T(D_M+\Theta_y)^{-1}\mathbf{w}}{\mathbf{w}^T(M+\Theta_y)^{-1}\mathbf{w}}\geq\lambda_{\min}\Big((D_M+\Theta_y)^{-1}(M+\Theta_y)\Big)\geq\min\left\{\lambda_{\min}(D_M^{-1}M),1\right\},
\end{equation*}
and hence that
\begin{equation}\label{omegabound}
\ \frac{\boldsymbol\omega^T\boldsymbol\omega}{\boldsymbol\gamma^T\boldsymbol\gamma}\geq{}R_{\Theta}\cdot\min\left\{\lambda_{\min}(D_M^{-1}M),1\right\}.
\end{equation}

Finally, we observe that $R_{\Theta}=\lambda_{\min}(D_M^{-1}M)=1$ for lumped mass matrices, as $D_M=M$. Inserting \eqref{omegabound} into \eqref{Rbound} then gives the required result.~~$\Box$

\vspace{1em}

\begin{remark} For consistent mass matrices, the working above still holds, except $R_{\Theta}$ and $\lambda_{\min}(D_M^{-1}M)$ are not equal to $1$. Therefore, the bound reads
\begin{equation*}
\ \lambda(\widehat{S}^{-1}S)\geq\frac{1}{2}\cdot\min\Big\{\min\hspace{0.1em}R_{\Theta}\cdot\min\left\{\lambda_{\min}(D_M^{-1}M),1\right\},1\Big\},
\end{equation*}
and depends on the matrix $[D_M-(\Theta+D_M^{-1})^{-1}]^{-1}[M-(\Theta+M^{-1})^{-1}]$, which does not have uniformly bounded eigenvalues. This is, however, a weak bound, and in practice we find that the (smallest and largest) eigenvalues of the preconditioned Schur complement are moderate in size.

Furthermore, in numerical experiments, we find the vast majority of the eigenvalues of $\widehat{S}^{-1}S$ 
to be clustered in the interval $\left [\frac{1}{2},1 \right]$, particularly as the Interior Point method approaches convergence, for the following reasons. In \cite[Theorem 4.1]{PW11}, it is shown that
\begin{eqnarray}\label{half1}
\ \lambda\left(\left[\left(L+\frac{1}{\sqrt{\alpha}}M\right)M^{-1}\left(L+\frac{1}{\sqrt{\alpha}}M\right)^T\right]^{-1}\left[LM^{-1}L^T+\frac{1}{\alpha}M\right]\right)\in\left[\frac{1}{2},1\right],
\end{eqnarray}
for any (positive) value of $\alpha$, and any mesh-size, provided $L+L^T$ is positive semidefinite, which is the case for Poisson and convection--diffusion problems for instance. For the Schur complement \eqref{Schur} and Schur complement approximation \eqref{SchurApprox}, as the Interior Point method approaches convergence, two cases will arise: (i) some entries of $\Theta_w^{-1}+\Theta_v^{-1}$ will approach zero, whereupon substituting these values into \eqref{Schur} and \eqref{SchurApprox} gives that $S$ and $\widehat{S}$ are both approximately $L(M+\Theta_y^{-1})^{-1}L^T$, so the eigenvalues of $\widehat{S}^{-1}S$ should be roughly $1$; (ii) some entries of $\Theta_w^{-1}+\Theta_v^{-1}$ approach infinity (with many entries of $\Theta_y$ correspondingly approaching zero), so $S$ is approximately $LM^{-1}L^T+\frac{1}{\alpha}M$, with $\widehat{S}$ an approximation of $(L+\frac{1}{\sqrt{\alpha}}M)M^{-1}(L+\frac{1}{\sqrt{\alpha}}M)^T$, giving clustered eigenvalues as predicted by \eqref{half1}.
The numerical evidence of the described behavior, for consistent mass matrices, is shown in Figure \ref{eig}.

\vspace{1em}

\tikzexternaldisable
 \begin{figure}[htb]
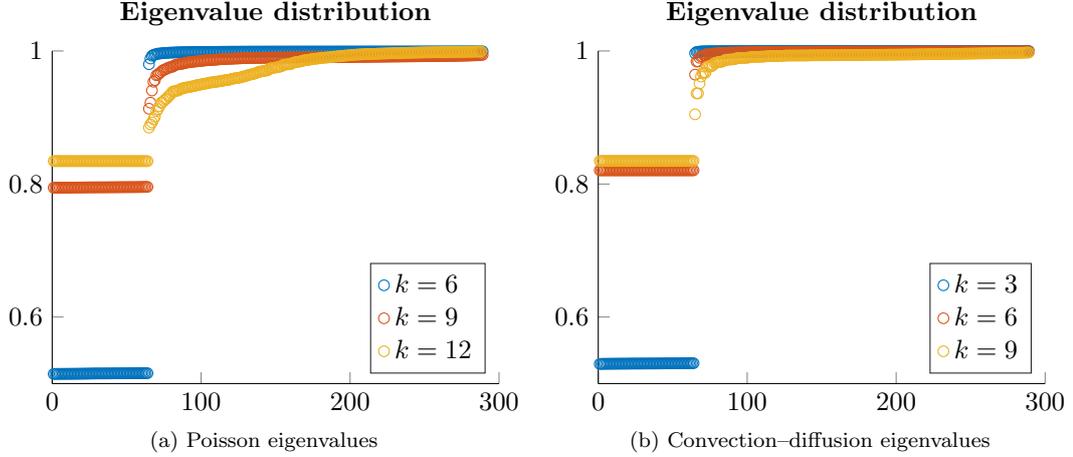

\begin{center}
	\setlength\figureheight{0.3\linewidth} 
	\setlength\figurewidth{0.4\linewidth}
	\subfloat[Poisson eigenvalues]{
	\input{figures/eigplot1.tikz}
	}
	\subfloat[Convection--diffusion eigenvalues]{
	\input{figures/eigplot2.tikz}
	}
 \end{center}
\caption{Eigenvalue distribution of $\widehat{S}^{-1}S $ at later Interior Point iterations for test problems involving Poisson's equation (left)
and the convection--diffusion equation (right) (with mesh-size $h=2^{-4}$). 
} 
\label{eig}
\end{figure}
\tikzexternalenable
\end{remark}

We note that the $(1,1)$-block and Schur complement approximations that we have derived are both symmetric positive definite, so we may apply the {\scshape minres} algorithm with a block diagonal preconditioner
of the form
\begin{equation*}
\ \mathcal{P}_D=\left[\begin{array}{cccc}
M+\Theta_y & 0 & 0 & 0 \\
0 & \alpha{}D_M+\Theta_w & -\alpha{}D_M & 0 \\
0 & -\alpha{}D_M & \alpha{}D_M+\Theta_v & 0 \\
0 & 0 & 0 & \widehat{S} \\
\end{array}\right],
\end{equation*}
with $\widehat{S}$ defined as above.

It is also possible to exploit the often faster convergence achieved by block triangular preconditioners within {\scshape gmres}, and utilize the block triangular preconditioner:
\begin{equation*}
\ \mathcal{P}_T=\left[\begin{array}{cccc}
M+\Theta_y & 0 & 0 & 0 \\
0 & \alpha{}D_M+\Theta_w & -\alpha{}D_M & 0 \\
0 & -\alpha{}D_M & \alpha{}D_M+\Theta_v & 0 \\
L & -M & M & -\widehat{S} \\
\end{array}\right].
\end{equation*}

\subsection{Preconditioner for Partial Observations}
\label{subsec::po}
In practice, the quantity of importance from a practical point-of-view is the difference between the state variable and the desired state on a certain region of the domain,
i.e. $\Omega_1\subset\Omega$, in which case one would instead consider the term $\frac{1}{2}\|\rm y-\rm y_d\|^ 2_{L^2(\Omega_1)}$ within the cost functional \eqref{pb}.
This results in a mass matrix where many of the eigenvalues are equal to zero. In more detail, the matrix $M+\Theta_y$ is in practice $M_s+\Theta_y$, where $M_s$ is a (singular) mass matrix acting on a subdomain, although for the purposes of our working we retain the existing notation. Hence, the standard saddle-point preconditioning
approach cannot be straightforwardly applied, due to the $(1,1)$-block being singular. One strategy is to replace the singular mass matrix with a slightly perturbed 
version in the preconditioning step. However, it is not straightforward to estimate the strength of this perturbation and its affect on the preconditioner.

Another alternative is presented in \cite{BenDOS15,herzog2018fast}, and we follow this strategy here. 
This method is tailored to the case where
the leading block of the saddle-point system is highly singular (meaning a large proportion of its eigenvalues are zero), due to the fact that the observations 
are placed only on parts of the domain.
In more detail, we consider the matrix system
\begin{equation}\label{MatrixPartial}
\left [\begin{array}{cc c}
 M+\Theta_y&0&L^T\\
 0&\alpha \widetilde M  + \Theta_z &-\bar M^T \\
 L&-\bar M&0\\         
\end{array} \right],
\end{equation}
with $M+\Theta_y$ often a highly singular matrix, as $\Theta_y=0$ when no state constraints are present. 
The mass matrix used to construct $\widetilde M$ is then defined on the control domain, which can be the whole domain or part of it. 
We start by considering the following permutation of the matrix to be solved:
\begin{equation}\label{Permuted}
	\Pi
\left [\begin{array}{ccc}
 M+\Theta_y&0&L^T\\
 0&\alpha \widetilde M  + \Theta_z &-\bar M^T \\
 L&-\bar M&0\\         
\end{array} \right]
=
\left [\begin{array}{ccc}
L&-\bar M&0\\     
0&\alpha \widetilde M  + \Theta_z &-\bar M^T \\
M+\Theta_y&0&L^T\\     
\end{array} \right]
,
\end{equation}
where 
\begin{equation*}
	\Pi:=
	\left[
		\begin{array}{ccc}
			0&0&I\\
			0&I&0\\
			I&0&0\\
		\end{array}
	\right].
\end{equation*}
The matrix \eqref{Permuted} is a block matrix of the form \eqref{ABCD} with
\begin{equation*}
	\ A=\left[
		\begin{array}{cc}
L&-\bar M\\     
0&\alpha \widetilde M  + \Theta_z\\
		\end{array}
	\right],\quad\quad{}B_1=\left[
		\begin{array}{cc}
			0\\
			-\bar M^T \\
		\end{array}
	\right],\quad\quad{}B_2=\left[
		\begin{array}{cc}
			M+\Theta_y&0\\
		\end{array}
	\right],\quad\quad{}C=\left[
		\begin{array}{c}
			L^T \\
		\end{array}
	\right],
\end{equation*}
which is a modification to a general saddle-point system, with non-symmetric extra-diagonal blocks and a non-zero $(2,2)$-block given by $L^T$. 
Based on this we propose the following preconditioner of block-triangular type for the permuted system:
\begin{equation*}
	\widetilde{\mathcal{P}}=
	\left[
		\begin{array}{ccc}
L&-\bar M&0\\     
0&\alpha \widetilde M  + \Theta_z &0\\
M+\Theta_y&0&-\widehat{S}_{\Pi}\\
		\end{array}
	\right],
\end{equation*}
with the inverse then given by
\begin{equation*}
	\widetilde{\mathcal{P}}^{-1}=
	\left[
		\begin{array}{ccc}
			L^{-1}&L^{-1}\bar M(\alpha \widetilde M  + \Theta_z )^{-1}&0\\
			0&(\alpha \widetilde M  + \Theta_z)^{-1} &0\\
			\widehat{S}_{\Pi}^{-1}(M+\Theta_y)L^{-1}&\widehat{S}_{\Pi}^{-1}(M+\Theta_y)L^{-1}\bar M(\alpha \widetilde M  + \Theta_z )^{-1}&-\widehat{S}_{\Pi}^{-1}\\
		\end{array}
	\right].
\end{equation*}
The matrix $\widehat{S}_{\Pi}$ is designed to approximate the Schur complement $S_{\Pi}$ of the \emph{permuted matrix system}, that is
\begin{equation*}
	\widehat{S}_{\Pi}\approx
	S_{\Pi}
	=L^{T}+(M+\Theta_y)L^{-1}\bar M(\alpha \widetilde M  + \Theta_z )^{-1}\bar M^T.
\end{equation*}
We now propose a preconditioner $\mathcal{P}_{\Pi}$ for the original matrix \eqref{MatrixPartial}, such that $\mathcal{P}_{\Pi}^{-1}=\widetilde{\mathcal{P}}^{-1}\Pi$, and we therefore obtain
\begin{equation}
	\label{eq:prec1}
	\mathcal{P}_{\Pi}^{-1}=
\left[
\begin{array}{ccc}
0&L^{-1}\bar M(\alpha \widetilde M  + \Theta_z )^{-1}&L^{-1}\\
0&(\alpha \widetilde M  + \Theta_z)^{-1} &0\\
-\widehat{S}_{\Pi}^{-1}&\widehat{S}_{\Pi}^{-1}(M+\Theta_y)L^{-1}\bar M(\alpha \widetilde M  + \Theta_z )^{-1}&\widehat{S}_{\Pi}^{-1}(M+\Theta_y)L^{-1}\\
\end{array}
\right].
\end{equation}
Applying the preconditioner is in fact more straightforward than it currently appears. To compute a vector $\mathbf{v}=\mathcal{P}_{\Pi}^{-1}\mathbf{w}$, where $\mathbf{v}:=\left[\mathbf{v}_{1}^T,~\mathbf{v}_{2}^T,~\mathbf{v}_{3}^T\right]^T$, $\mathbf{w}:=\left[\mathbf{w}_{1}^T,~\mathbf{w}_{2}^T,~\mathbf{w}_{3}^T\right]^T$, we first observe from the second block of $\mathcal{P}_{\Pi}^{-1}$ that
\begin{equation*}
	(\alpha \widetilde M  + \Theta_z)^{-1}\mathbf{w}_2=\mathbf{v}_2.
\end{equation*}
The first equation derived from \eqref{eq:prec1} then gives that
\begin{align*}
L^{-1}\bar M(\alpha \widetilde M  + \Theta_z )^{-1}\mathbf{w}_2+L^{-1}\mathbf{w}_3&=\mathbf{v}_1\\
\Rightarrow\hspace{7.2em}L^{-1}(\bar M\mathbf{v}_2+\mathbf{w}_3)&=\mathbf{v}_1,
\end{align*}
and applying this within the last equation in \eqref{eq:prec1} that
\begin{align*}
-\widehat{S}_{\Pi}^{-1}\mathbf{w}_1+\widehat{S}_{\Pi}^{-1}(M+\Theta_y)L^{-1}\bar M(\alpha \widetilde M  + \Theta_z )^{-1}\mathbf{w}_2+\widehat{S}_{\Pi}^{-1}(M+\Theta_y)L^{-1}\mathbf{w}_3&=\mathbf{v}_3\\
\Rightarrow\hspace{5.1em}-\widehat{S}_{\Pi}^{-1}\mathbf{w}_1+\widehat{S}_{\Pi}^{-1}(M+\Theta_y)\big(L^{-1}\bar M(\alpha \widetilde M  + \Theta_z )^{-1}\mathbf{w}_2+L^{-1}\mathbf{w}_3\big)&=\mathbf{v}_3\\
\Rightarrow\hspace{20.95em}\widehat{S}_{\Pi}^{-1}\big((M+\Theta_y)\mathbf{v}_1-\mathbf{w}_1\big)&=\mathbf{v}_3.
\end{align*}

Thus we need to approximately solve with $\widehat{S}_{\Pi}$, $L$, and $\alpha \widetilde M  + \Theta_z$, which are all invertible matrices, to apply the preconditioner. We now briefly discuss our choice of $\widehat{S}_{\Pi}.$ We suggest a matching strategy as above, to write
\begin{align*}
S_{\Pi}=L^{T}+(M+\Theta_y)L^{-1}\bar M(\alpha \widetilde M  + \Theta_z )^{-1}\bar M^T\approx\big(L^{T}+{M}_l\big)L^{-1}\big(L+{M}_r\big)=\widehat{S}_{\Pi},
\end{align*}
where 
\begin{equation*}
	{M}_lL^{-1}{M}_r\approx(M+\Theta_y)L^{-1}\bar M(\alpha \widetilde M  + \Theta_z )^{-1}\bar M^T.
\end{equation*}
Such an approximation may be achieved if, for example, 
\begin{equation*}
	{M}_l=M+\Theta_y,\quad\quad{M}_r\approx \bar M(\alpha \widetilde M  + \Theta_z )^{-1}\bar M^T.
\end{equation*}
Alternatively, we can use a matrix based on the approximation $\widehat{M}$ from the previous section to approximate ${M}_r.$
We thus build such approximations into our preconditioner $\mathcal{P}_{\Pi}$, although further tailoring of such preconditioners is a subject of future investigation.


\subsection{Time-Dependent Problems}
\label{subsec::td}
To demonstrate the applicability of our preconditioners to time-dependent PDE-constrained optimization problems, we now consider the minimization of the cost functional
\begin{equation*}
\ \F(\rm y,\rm u)=\frac{1}{2}\|\rm y-\rm y_d\|^ 2_{L^2(\Omega\times(0,T))}+ \frac{\alpha}{2}\|\rm u\|^ 2_{L^2(\Omega\times(0,T))} + \beta\|u\|_{L^1(\Omega\times(0,T))},
\end{equation*}
subject to the PDE $\rm y_{t}-\Delta\rm y=\rm u+\rm f$ on the space-time interval $\Omega\times(0,T)$, along with suitable boundary and initial conditions.

With the backward Euler method used to handle the time derivative, the matrix within the system to be solved is of the form
\begin{equation}\label{TimeDeptSystem}
 \mathcal{A} = \left [\begin{array}{c c  c }

 \tau \mathcal{M}_c   + \Theta_y         &  0               &   \mathcal{L}^T    \\
 0             &  \alpha\tau\widetilde{\mathcal{M}}_c  + \Theta_z &   -\tau\bar{\mathcal{M}}^T     \\
 \mathcal{L}             &  -\tau\bar{\mathcal{M}}             &    0   \\         
        \end{array} \right],
\end{equation}
with $\tau$ the time-step used.

The matrix $\mathcal{M}_c$ is a block diagonal matrix consisting of multiples of mass matrices on each block diagonal corresponding to each time-step, depending on the quadrature rule used to approximate the cost functional in the time domain. For example, if a trapezoidal rule is used, then $\mathcal{M}_c=\text{blkdiag}(\frac{1}{2}M,M,...,M,\frac{1}{2}M)$, and if a rectangle rule is  used, then $\mathcal{M}_c=\mathcal{M}:=\text{blkdiag}(M,M,...,M,M)$. Further,
\begin{equation*}
\ \widetilde{\mathcal{M}}_c=\left[\begin{array}{cc}
\mathcal{M}_c & -\mathcal{M}_c \\
-\mathcal{M}_c & \mathcal{M}_c \\
\end{array}\right],\quad\quad\bar{\mathcal{M}}=\left[\begin{array}{cc}
\mathcal{M} & -\mathcal{M} \\
\end{array}\right],
\end{equation*}
and $\mathcal{L}$ is defined as follows (with its dimension equal to that of $L$, multiplied by the number of time-steps):
\begin{equation*}
\ \mathcal{L}=\left[\begin{array}{cccc}
M+\tau{}L & & & \\
-M & M+\tau{}L & & \\
 & \ddots & \ddots & \\
 & & -M & M+\tau{}L \\
\end{array}\right].
\end{equation*}

We now consider saddle-point preconditioners for the matrix \eqref{TimeDeptSystem}. We may apply preconditioners of the form
\begin{align*}
\ \mathcal{P}_{D}={}&\left [\begin{array}{c c c c}
\tau \mathcal{M}_c + \Theta_y  & 0 & 0 & 0 \\
0 & \alpha\tau\mathcal{D}_{M_c}  + \Theta_w & -\alpha\tau\mathcal{D}_{M_c} & 0 \\
0 & -\alpha\tau\mathcal{D}_{M_c} & \alpha\tau\mathcal{D}_{M_c}  + \Theta_v & 0 \\
0 & 0 & 0 & \widehat{\mathcal{S}} \\         
\end{array} \right] \\
\ \text{or}\quad\mathcal{P}_{T}={}&\left [\begin{array}{c c c c}
\tau \mathcal{M}_c + \Theta_y  & 0 & 0 & 0 \\
0 & \alpha\tau\mathcal{D}_{M_c}  + \Theta_w & -\alpha\tau\mathcal{D}_{M_c} & 0 \\
0 & -\alpha\tau\mathcal{D}_{M_c} & \alpha\tau\mathcal{D}_{M_c}  + \Theta_v & 0 \\
\mathcal{L} & -\tau\mathcal{M} & \tau\mathcal{M} & -\widehat{\mathcal{S}} \\         
\end{array} \right],
\end{align*}
where $\mathcal{D}_{M_c}:=\text{diag}(\mathcal{M}_c)$, the matrix $\tau\mathcal{M}_{c}+\Theta_y$ can be approximately inverted by applying Chebyshev semi-iteration to the matrices arising at each time-step, and $\widehat{\mathcal{S}}$ is an approximation of the Schur complement:
\begin{equation*}
\ \mathcal{S}=\mathcal{L}(\tau\mathcal{M}_{c}+\Theta_y)^{-1}\mathcal{L}^{T}+\frac{\tau}{\alpha}\mathcal{M}\mathcal{M}_c^{-1}\mathcal{M}-\frac{1}{\alpha^2}\mathcal{M}\mathcal{M}_c^{-1}\left(\Theta_w^{-1}+\Theta_v^{-1}+\frac{1}{\alpha\tau}\mathcal{M}_c^{-1}\right)\mathcal{M}_c^{-1}\mathcal{M}.
\end{equation*}
We select the approximation
\begin{equation*}
\ \widehat{\mathcal{S}}=\left(\mathcal{L}+\widehat{\mathcal{M}}\right)(\tau\mathcal{M}_{c}+\Theta_y)^{-1}\left(\mathcal{L}+\widehat{\mathcal{M}}\right)^{T},
\end{equation*}
using the same reasoning as in Section \ref{sec:Schur}, where
\begin{equation*}
\ \widehat{\mathcal{M}}=\left[\frac{\tau}{\alpha}\mathcal{D}_{M}^2\mathcal{D}_{M_c}^{-1}-\frac{1}{\alpha^2}\mathcal{D}_{M}^2\mathcal{D}_{M_c}^{-2}\left(\Theta_w^{-1}+\Theta_v^{-1}+\frac{1}{\alpha\tau}\mathcal{D}_{M_c}^{-1}\right)\right]^{1/2}(\tau\mathcal{D}_{M_c}+\Theta_y)^{1/2},
\end{equation*}
with $\mathcal{D}_{M}:=\text{diag}(\mathcal{M})$. Within the numerical experiments of the forthcoming section, we apply the preconditioning strategy that arises from the working above.

\section{Numerical Experiments}\label{exp}

We now implement the Interior Point algorithm described in the Appendix, using {\scshape matlab}\textsuperscript{\textregistered} R2017b 
on an Intel\textsuperscript{\textregistered} Xeon\textsuperscript{\textregistered} computer with a 2.40GHz processor, and 250GB of RAM.
Within the algorithm we employ the preconditioned {\scshape minres}\ \cite{minres} and {\scshape gmres} \cite{gmres} methods with the following preconditioners:
\begin{itemize}
\item \ipmbt: {\sc gmres} and block triangular preconditioner $\mathcal{P}_T,$ 
\item \ipmbd: {\sc minres} with block diagonal preconditioner  $\mathcal{P}_D,$
\item \ipmbpi: {\sc gmres} and block triangular preconditioner $\mathcal{P}_\Pi.$
\end{itemize}
Regarding the parameters listed in the Appendix, we use
$\alpha_0 = 0.995$ and $\epsilon_p=\epsilon_d=\epsilon_c = 10^{-6}$.
For the barrier reduction parameter $\sigma$, we consider for each class of
problems tested a value that ensures a smooth decrease in the complementarity measure
$\xi^k_c$ in (\ref{gap}), that is to say $\|\xi^k_c\| = \mathcal{O}(\mu^k)$. This way, the number of 
nonlinear (Interior Point) iterations typically depends only on $\sigma$.
We solve the linear matrix systems to a (relative unpreconditioned residual norm) tolerance of $10^{-10}$.

\begin{figure}
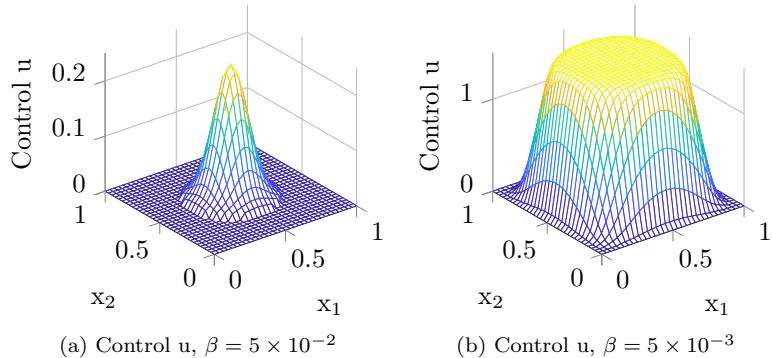

\begin{center}
	\setlength\figureheight{0.225\linewidth} 
	\setlength\figurewidth{0.225\linewidth} 
 \subfloat[Control $\rm u$, $\beta=5\times10^{-2}$]{
	\input{figures/controlPoissonbeta5e_2.tikz}
	}
	\subfloat[Control $\rm u$, $\beta=5\times10^{-3}$]{
	\input{figures/controlPoissonbeta5e_3.tikz}
	}
 \end{center}
\caption{Poisson problem: computed solutions of the control $\rm u$, for two values of $\beta$.} \label{fig::poissonu}
\end{figure}

\begin{table}[htb!]
\begin{center}
\begin{tabular}{ccccccc}
\toprule
           & \multicolumn{ 2}{c}{$\beta = 10^{-1}$} & \multicolumn{ 2}{c}{$\beta =   10^{-2}$} & \multicolumn{ 2}{c}{$\beta = 10^{-3}$} \\
\midrule
           &   {\sc sparsity} &  $\|u\|_1$ &   {\sc sparsity} &  $\|u\|_1$ &  {\sc sparsity} &  $\|u\|_1$ \\
\midrule
$\alpha = 10^{-2}$ &      99\%  &       3 &       15\% &       $7\times 10^2$ &      12\%  &       $1\times 10^3$ \\

$\alpha = 10^{-4}$ &      100\% &       2 &       38\% &       $9\times 10^2$ &      12\%  &       $1\times 10^3$ \\

$\alpha = 10^{-6}$ &      100\% &       2 &       39\% &       $9\times 10^2$ &      12\%  &       $1\times 10^3$ \\

\bottomrule
\end{tabular}  
\end{center}
\caption{Poisson problem: sparsity features of the computed optimal control, for a range of $\alpha$ and $\beta$, and mesh-size $h = 2^{-5}$. 
\label{tab::sparsity}}
\end{table}

We apply the {\scshape ifiss} software package \cite{ifissmatlab,ifisslink} to build
the relevant finite element matrices for the 2D examples shown in this section, and use the
{\scshape deal.II} library \cite{dealii} in the 3D case. In each case we utilize $Q1$ finite elements
for the state, control, and adjoint variables.

We apply $20$ steps of Chebyshev semi-iteration to approximate the inverse of mass matrices, as well as mass matrices plus positive diagonal matrices, whenever they arise within the preconditioners.
Applying the approximate inverses of the Schur complement approximations derived for each of our preconditioners
requires solving for matrices of the form $L + \widehat M$ and its transpose.
For this we utilize $3$ V-cycles of the algebraic multigrid routine {\sc hsl-mi20} \cite{Boyle2007},
with a Gauss--Seidel coarse solver, and apply $5$ steps of pre- and post-smoothing. 
For time-dependent problems, we also use Chebyshev semi-iteration and algebraic multigrid within the preconditioner, 
but are required to apply the methods to matrices arising from each time-step.
In all the forthcoming tables of results, we report the average number of linear ({\scshape minres} or {\scshape gmres}) iterations {\sc av-li},
and the average CPU time {\sc av-cpu}. The overall number of nonlinear (Interior Point) iterations {\sc nli} is specified in the table captions. 
We believe these demonstrate the effectiveness of our proposed Interior Point and preconditioning approaches, as well as the robustness of the
overall method, for a range of PDEs, matrix dimensions, and parameters involved in the problem set-up.

\subsection{A Poisson Problem}

\begin{table}[htb!]
\begin{center}
\begin{tabular}{llrrrrrr}
\toprule
             &                               &      \multicolumn{ 2}{c}{\ipmbt }    &     \multicolumn{ 2}{|c }{\ipmbd}             \\
\midrule
       $h=2^{-\ell}$ &   $\mathrm{log}_{10}\alpha$ &  {\sc av-li}  &  {\sc av-cpu} &      {\sc av-li}  &  {\sc av-cpu}  \\
\midrule

\multicolumn{ 1}{c}{6} &         $-2$ &        8.9 &        0.2 &       19.4 &        0.4 \\

\multicolumn{ 1}{c}{} &         $-4$ &        7.2 &        0.2 &       16.3 &        0.3 \\

\multicolumn{ 1}{c}{} &         $-6$ &        7.1 &        0.2 &       14.6 &        0.3 \\

\multicolumn{ 1}{c}{7} &         $-2$ &        9.0 &        0.8 &       19.5 &        1.6 \\

\multicolumn{ 1}{c}{} &         $-4$ &        7.1 &        0.7 &       15.8 &        1.3 \\

\multicolumn{ 1}{c}{} &         $-6$ &        6.8 &        0.6 &       14.4 &        1.4 \\

\multicolumn{ 1}{c}{8} &         $-2$ &        6.9 &        2.5 &       14.3 &        5.0 \\

\multicolumn{ 1}{c}{} &         $-4$ &        6.5 &        2.4 &       13.4 &        4.7 \\

\multicolumn{ 1}{c}{} &         $-6$ &        6.5 &        2.4 &       12.8 &        4.5 \\

\multicolumn{ 1}{c}{9} &         $-2$ &        7.9 &       12.4 &       13.8 &       21.8 \\

\multicolumn{ 1}{c}{} &         $-4$ &        7.6 &       12.0 &       12.7 &       20.2 \\

\multicolumn{ 1}{c}{} &         $-6$ &        7.5 &       11.9 &       12.3 &       20.0 \\
\bottomrule
\end{tabular}  
\end{center}
\caption{Poisson problem: average Krylov iterations and CPU times for problem with control constraints, for a range of $h$ and $\alpha$, $\beta = 10^{-2}$, $\sigma = 0.2$, $\textsc{nli} = 9$.
\label{tab::resultspoisson1}}
\end{table}

We first examine an optimization problem involving Poisson's equation, investigating the behavior of the IPM and our proposed preconditioners. 

%

\subsection*{Two-Dimensional Case}
We focus initially on the performance of our solvers for the two-dimensional Poisson problem, employing both \ipmbt and \ipmbd methods, as well as considering some sparsity issues. 
We set the box constraints for the control to be $\rm u_a=-2$, $\rm u_b=1.5,$ and the desired state 
$\rm y_d=\sin(\pi {\rm x_1})\sin(\pi {\rm x_2}) $, with ${\rm x}_i$ denoting the $i$th spatial variable. Figure \ref{fig::poissonu} displays the computed optimal controls for this problem for a particular set-up on the domain  $\Omega=(0,1)^2$, for both $\beta=5\times10^{-2}$ and $\beta=5\times10^{-3}$
as well as $\alpha = 10^{-2}$. Table \ref{tab::sparsity} reports the level of sparsity in the computed solution, as well as its 
$\ell_1$-norm, when varying the regularization parameters $\alpha$ and $\beta$. The value of {\sc sparsity} in the table is computed by
measuring the percentage of components of $u$ which are below a certain threshold ($10^{-2}$ in our case),
see e.g. \cite{fpcas}. We observe that our algorithm reliably computes sparse
controls, and as expected the sparsity of the solution increases when $\beta$ is correspondingly increased.

\begin{figure}[htb]
   \centering
	\setlength\figureheight{0.35\linewidth} 
	\setlength\figurewidth{0.45\linewidth}
%
%
\definecolor{mycolor1}{rgb}{0.00000,0.44700,0.74100}%
\definecolor{mycolor2}{rgb}{0.85000,0.32500,0.09800}%
\definecolor{mycolor3}{rgb}{0.92900,0.69400,0.12500}%
\definecolor{mycolor4}{rgb}{0.49400,0.18400,0.55600}%
\begin{tikzpicture}

\begin{axis}[%
width=\figurewidth,
height=\figureheight,
at={(0.747in,0.477in)},
scale only axis,
xmin=0,
xmax=12,
xlabel={$k$},
ymode=log,
ymin=1e-15,
ymax=100000,
yminorticks=true,
axis background/.style={fill=white},
legend style={at={(0.732,0.62)},anchor=south west,legend cell align=left,align=left,draw=white!15!black}
]
\addplot [color=mycolor1,solid,mark=asterisk,mark options={solid}]
  table[row sep=crcr]{%
1	1\\
2	0.25\\
3	0.0625\\
4	0.01562\\
5	0.003906\\
6	0.0009766\\
7	0.0002441\\
8	6.104e-05\\
9	1.526e-05\\
10	3.815e-06\\
11	9.537e-07\\
12	2.384e-07\\
};
\addlegendentry{$\mu{}^{k}$};

\addplot [color=mycolor2,solid,mark=asterisk,mark options={solid}]
  table[row sep=crcr]{%
1	16.12\\
2	0.0806\\
3	0.000403\\
4	2.015e-06\\
5	1.007e-08\\
6	5.04e-11\\
7	2.523e-13\\
8	9.913e-14\\
9	1.046e-14\\
10	1.051e-14\\
11	5.065e-15\\
12	3.539e-15\\
};
\addlegendentry{$\xi{}^{k}_{p}$};

\addplot [color=mycolor3,solid,mark=asterisk,mark options={solid}]
  table[row sep=crcr]{%
1	15.74\\
2	0.07873\\
3	0.0003935\\
4	1.961e-06\\
5	2.685e-08\\
6	6.398e-09\\
7	1.619e-09\\
8	4.076e-10\\
9	1.015e-10\\
10	2.44e-11\\
11	5.108e-12\\
12	5.787e-13\\
};
\addlegendentry{$\xi{}^{k}_{d}$};

\addplot [color=mycolor4,solid,mark=asterisk,mark options={solid}]
  table[row sep=crcr]{%
1	331.6\\
2	9.109\\
3	0.2646\\
4	0.03145\\
5	0.007769\\
6	0.001942\\
7	0.0004855\\
8	0.0001214\\
9	3.036e-05\\
10	7.602e-06\\
11	1.913e-06\\
12	4.896e-07\\
};
\addlegendentry{$\xi{}^{k}_{c}$};

\end{axis}
\end{tikzpicture}%
   \caption{Typical convergence history of the relevant quantities $\mu^k, \xi^k_p, \xi^k_d, \xi^k_c$. 
\label{fig::convh}}
\end{figure}
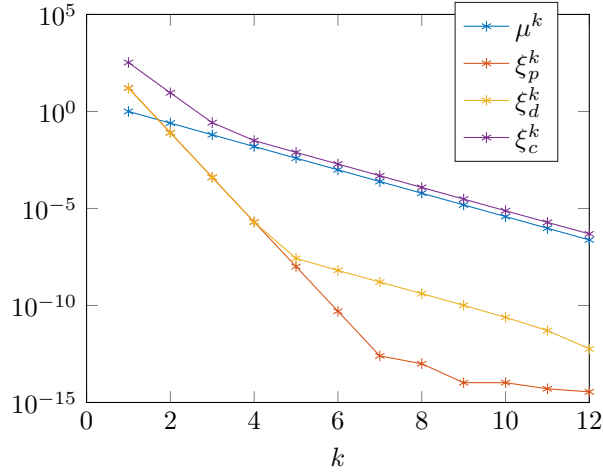

In Table \ref{tab::resultspoisson1} we compare the performance of the preconditioners $\mathcal{P}_T$ and $\mathcal{P}_D$ within the IPM, varying the 
spatial mesh-size $h = 2^{-i},\ i = 6, \dots, 9$, as well as the regularization parameter $\alpha$, while fixing the value $\beta = 10^{-2}$ (Table \ref{tab::sparsity} indicates that this value of $\beta$ gives rise to the most computationally interesting case). We set $\sigma = 0.2$, and
take $9$ Interior Point iterations with a final value $\mu^k = 5 \times 10^{-7}$. Figure \ref{fig::convh} provides a representation of the typical convergence behavior for the feasibilities $\xi^k_p, \xi^k_d$ and complementarity $\xi^k_c$, together with 
the decrease of $\mu^k$ with this value of $\sigma$.
The reported results demonstrate good robustness of both preconditioners with respect to both $h$ and $\alpha$ in terms of linear iterations and
CPU time, with \ipmbt outperforming \ipmbd in each measure.
Despite the fact that the value of {\sc av-li} is constant in both implementations, we observe that when using \ipmbd the number of
preconditioned {\scshape minres} iterations slightly increases as $\mu^k \rightarrow 0$, as many entries of $\Theta_{z}$ tend to zero. 
On the contrary, the number of preconditioned {\scshape gmres} iterations hardly varies with $k$.

\tikzexternaldisable
\tikzexternalenable

As a final validation of the general framework outlined, we report in Table \ref{tab::resultspoisson2}
results obtained when imposing both control and state constraints within the Poisson setting described above.
In particular, we set  $\rm y_a=-0.1$, $\rm y_b=0.8$, $\rm u_a=-1$, $\rm u_b=15$ and test the most promising implementation
of the IPM, that is the \ipmbt routine, while varying $h$ and $\alpha$. The reported values of {\sc av-li} confirm the roboustness of
the preconditioning strategy proposed.

\begin{table}[htb!]
\begin{center}
\begin{tabular}{llrr}
\toprule
             &                                   &             \multicolumn{ 2}{c}{\ipmbt} \\ 
\midrule
       $h=2^{-\ell}$ &   $\mathrm{log}_{10}\alpha$ &   {\sc av-li}  &  {\sc av-cpu} \\
\midrule
\multicolumn{ 1}{c}{6} &         $-2$ &       15.8 &        0.4 \\

\multicolumn{ 1}{c}{} &         $-4$ &       11.4 &        0.3 \\

\multicolumn{ 1}{c}{} &         $-6$ &       10.6 &        0.2 \\

\multicolumn{ 1}{c}{7} &         $-2$ &       14.8 &        1.5 \\

\multicolumn{ 1}{c}{} &         $-4$ &       11.4 &        1.0 \\

\multicolumn{ 1}{c}{} &         $-6$ &       10.3 &        0.9 \\

\multicolumn{ 1}{c}{8} &         $-2$ &       14.6 &        5.4 \\

\multicolumn{ 1}{c}{} &         $-4$ &       10.8 &        3.9 \\

\multicolumn{ 1}{c}{} &         $-6$ &       10.1 &        3.5 \\

\multicolumn{ 1}{c}{9} &         $-2$ &       14.5 &       22.1 \\

\multicolumn{ 1}{c}{} &         $-4$ &       10.8 &       16.6 \\

\multicolumn{ 1}{c}{} &         $-6$ &        9.0 &       15.4 \\

\bottomrule
\end{tabular}  \hfill \begin{tabular}{clrrrr}
\toprule
&&\multicolumn{ 2}{c}{\ipmbpi}\\
\midrule
       $h=2^{-\ell}$ &   $\mathrm{log}_{10}\alpha$ &  {\sc av-li}  &  {\sc av-cpu}  \\
\midrule

\multicolumn{ 1}{c}{3} &         $-2$ &       10.2 &       0.04 \\

\multicolumn{ 1}{c}{} &         $-4$ &       11.3 &       0.05 \\

\multicolumn{ 1}{c}{} &         $-6$ &       11.3 &       0.05 \\

\multicolumn{ 1}{c}{4} &         $-2$ &       11.2 &        0.4 \\

\multicolumn{ 1}{c}{} &         $-4$ &       11.3 &        0.4 \\

\multicolumn{ 1}{c}{} &         $-6$ &       11.3 &        0.4 \\

\multicolumn{ 1}{c}{5} &         $-2$ &       15.0 &        7.2 \\

\multicolumn{ 1}{c}{} &         $-4$ &       15.1 &        7.3 \\

\multicolumn{ 1}{c}{} &         $-6$ &       15.1 &        7.3 \\

\bottomrule
\end{tabular}  
\end{center}
\caption{\emph{(Left)} Poisson problem: average Krylov iterations and CPU times for problem with both control and state constraints, for a range of $h$ and $\alpha$, $\beta = 10^{-2}$, $\sigma = 0.2$ ($\textsc{nli} = 14$).\\\emph{(Right)} Three-dimensional Poisson problem with partial observations: average Krylov iterations and CPU times for problem, for a range of $h$ and $\alpha$, $\beta = 10^{-3}$, $\sigma = 0.25$ ($\textsc{nli} = 11$).
\label{tab::resultspoisson2}}
\end{table}

\subsection*{Three-Dimensional Case with Partial Observations}
\begin{figure}[htb!]
\begin{center}
 \setlength\figureheight{0.33\linewidth} 	
	\subfloat[Computed control $\rm u$]{
 		\includegraphics[width=0.33\textwidth]{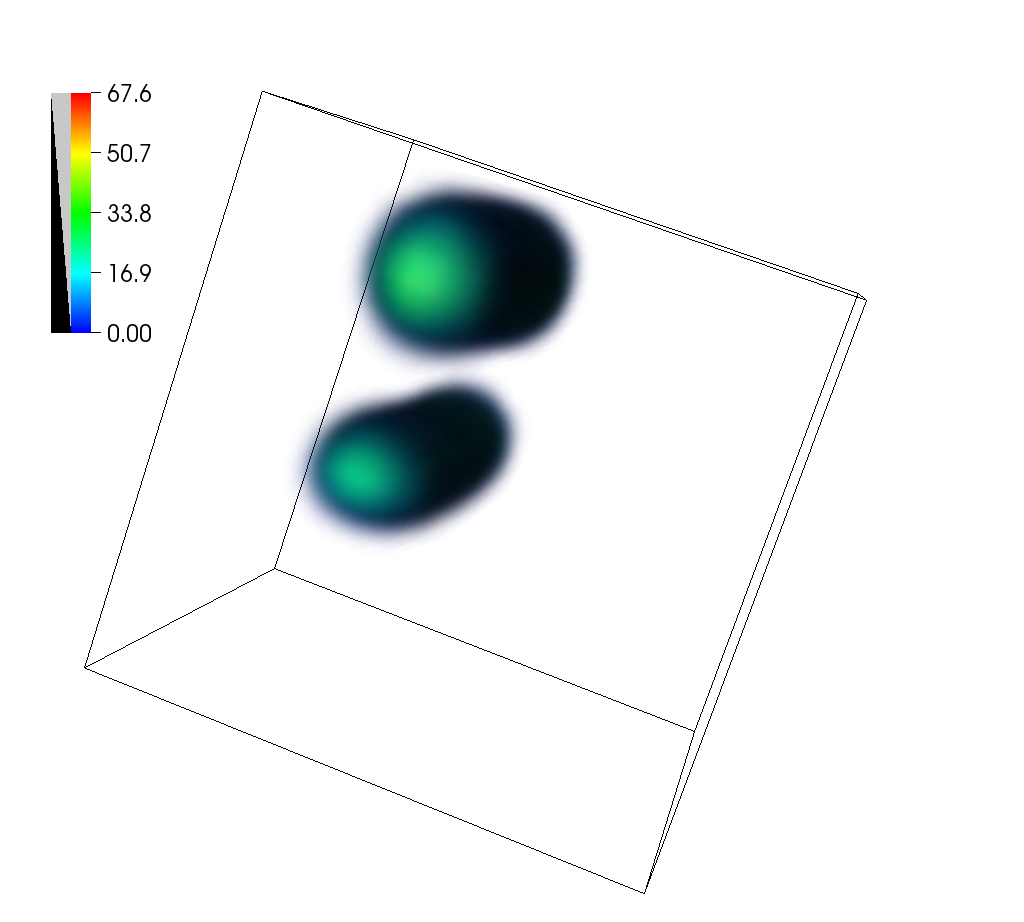}
	}
	\subfloat[Computed state $\rm y$]{
		\includegraphics[width=0.33\textwidth]{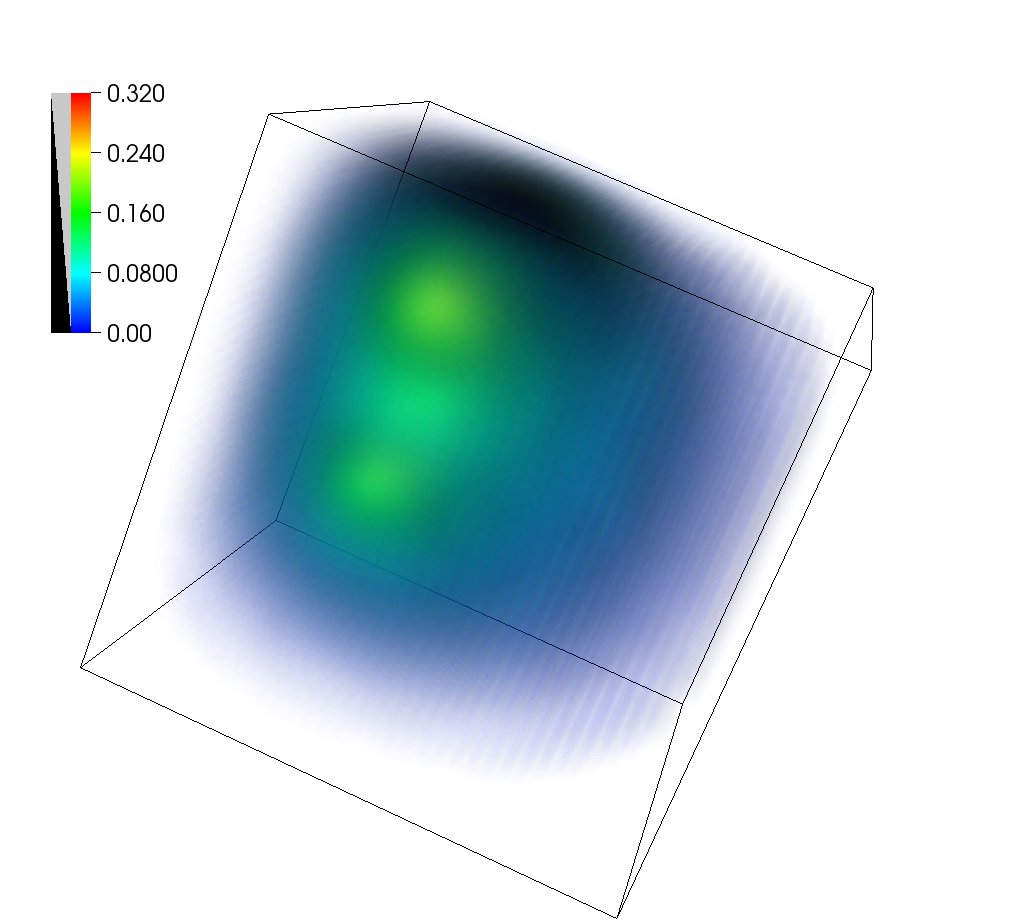}
	}
	\subfloat[Desired state $\rm y_d$]{
		\includegraphics[width=0.33\textwidth]{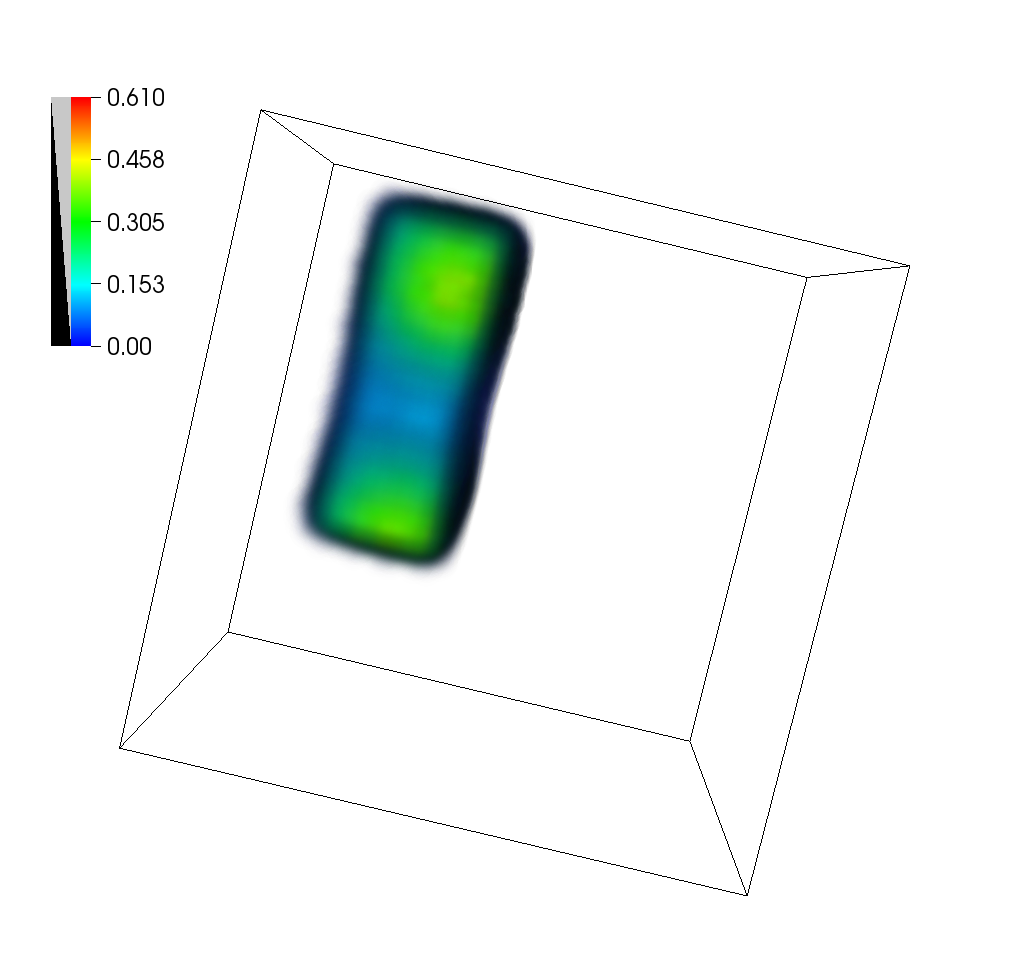}
	}
\end{center}

	\caption{Three-dimensional Poisson problem with partial observations: computed solutions for the control, state, and desired state.}\label{fig::res3d}
\end{figure}

We also wish to present results for the case of partial observations, paired with a three-dimensional example involving Poisson's equation on $\Omega = (0,1)^3$. 
The desired state is illustrated in Figure \ref{fig::res3d}. We use the preconditioner $\ppi$, as the observation domain {\color{black} $\Omega_1$}  is given by $0.2<\rm x_1<0.4$, $0.4<\rm x_2<0.9$, $0\leq\rm x_3\leq 1$,
and therefore the $(1,1)$-block of the matrix \eqref{NewtonSystem} is singular. The results for the computation with $\alpha = 10^{-5},$ $\beta = 10^{-3},$ and without additional box constraints, are also presented in Figure \ref{fig::res3d}, with the discretization involving $35937$ degrees of freedom.
To illustrate the performance of the proposed preconditioner $\mathcal{P}_\Pi$ with respect to changes in the parameter regimes, in Table \ref{tab::resultspoisson2} we provide results for a computation involving sparsity constraints applied to the control, as well as partial observation of the state, and set $\rm u_a=-2$, $\rm u_b=1.5.$ 
Again, the results are very promising and a large degree of robustness is achieved.

\subsection{A Convection--Diffusion Problem}
We next consider the optimal control of the convection--diffusion equation given by
$- \eps \Delta {\rm y} + \vec{\rm w} \cdot \nabla {\rm y} = {\rm u}$
on the domain $\Omega=(0,1)^2$, with the wind vector $\vec{\rm w}$ given by $\vec{\rm w} = \big[{\rm 2x_2(1-x_1^2)}, {\rm -2x_1(1-x_2^2)}\big]^T$, and the bounds on the control given by $\rm u_a=-2$ and $\rm u_b = 1.5$.
The desired state is here defined by
$\rm y_d = \exp(\rm -64(x_1-0.5)^2+(x_2-0.5)^2)$.
The discretization is again performed using Q1 finite elements, while also employing the Streamline Upwind Petrov--Galerkin (SUPG) \cite{BroH82} upwinding scheme as implemented in {\scshape ifiss}. The results of our scheme are given in Table \ref{tab::resultscd1}, which again exhibit robustness with respect to $h$ and $\alpha$, while also performing well for both values of $\eps$ tested.
\begin{figure}
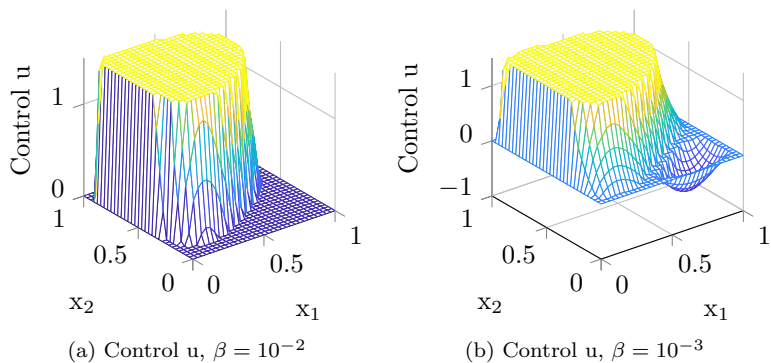

\begin{center}
	\setlength\figureheight{0.225\linewidth} 
	\setlength\figurewidth{0.225\linewidth} 

		\subfloat[Control $\rm u$, $\beta=10^{-2}$]{
	\input{figures/controlCDbeta1e_2.tikz}
	}
	\subfloat[Control $\rm u$, $\beta=10^{-3}$]{
	\input{figures/controlCDbeta1e_3.tikz}
	}
 \end{center}
\caption{Convection--diffusion problem: computed solutions of the control $\rm u$, for two values of $\beta$.} \label{fig::CDu}
\end{figure}

\begin{table}[htb!]
\begin{center}
\begin{tabular}{llrrrr|rrrr}
\toprule
           &            & \multicolumn{ 4}{c|}{$\eps = 10^{-1}$ } & \multicolumn{ 4}{c}{$\eps = 10^{-2}$} \\
					\midrule
             &                               &      \multicolumn{ 2}{c}{\ipmbt}    &     \multicolumn{ 2}{|c }{\ipmbd}
												 &      \multicolumn{ 2}{|c}{\ipmbt}    &     \multicolumn{ 2}{|c }{\ipmbd}\\
\midrule
       $h=2^{-\ell}$ &   $\mathrm{log}_{10}\alpha$ &  {\sc av-li}  &  {\sc av-cpu} &      {\sc av-li}  &  {\sc av-cpu} &  {\sc av-li}  &  {\sc av-cpu} &      {\sc av-li}  &  {\sc av-cpu} \\
\midrule
\multicolumn{ 1}{c}{6} &         $-2$ &       9.4  &       0.2  &       21.1 &        0.5 &       11.2 &        0.5 &       25.8 &        1.1 \\

\multicolumn{ 1}{c}{} &         $-4$ &        8.3 &        0.2 &       18.2 &        0.4 &       10.5 &        0.5 &       23.2 &        1.0 \\

\multicolumn{ 1}{c}{} &         $-6$ &        8.2 &        0.2 &       17.8 &        0.4 &       10.5 &        0.5 &       23.5 &        1.0 \\

\multicolumn{ 1}{c}{7} &         $-2$ &        8.2 &        0.8 &       18.0 &        1.7 &        9.2 &        1.6 &       20.6 &        3.4 \\

\multicolumn{ 1}{c}{} &         $-4$ &        7.5 &        0.7 &       16.3 &        1.5 &        8.7 &        1.5 &       19.0 &        3.1 \\

\multicolumn{ 1}{c}{} &         $-6$ &        7.5 &        0.7 &       16.1 &        1.5 &        8.7 &        1.5 &       19.4 &        3.1 \\

\multicolumn{ 1}{c}{8} &         $-2$ &        7.5 &        2.7 &       16.3 &        5.6 &        8.0 &        3.8 &       17.1 &        7.9 \\

\multicolumn{ 1}{c}{} &         $-4$ &        7.0 &        2.5 &       15.1 &        5.2 &        7.7 &        3.7 &       16.4 &        7.5 \\

\multicolumn{ 1}{c}{} &         $-6$ &        7.0 &        2.5 &       14.8 &        5.1 &        7.7 &        3.7 &       16.4 &        7.5 \\

\multicolumn{ 1}{c}{9} &         $-2$ &        7.0 &       11.2 &       14.9 &       23.0 &        7.3 &       13.1 &       15.1 &       26.3 \\

\multicolumn{ 1}{c}{} &         $-4$ &        6.7 &       11.0 &       14.2 &       22.4 &        6.8 &       12.5 &       14.4 &       25.5 \\

\multicolumn{ 1}{c}{} &         $-6$ &        6.7 &       11.0 &       13.9 &       21.7 &        6.8 &       12.5 &       14.5 &       25.5 \\

\bottomrule
\end{tabular}  
\end{center}
\caption{Convection--diffusion problem: average Krylov iterations and CPU times for problem with control constraints, for a range of $h$ and $\alpha$, $\beta = 10^{-3}$, $\sigma=0.25$ ($\textsc{nli} = 11$) with $\eps = 10^{-1}$, and $\sigma=0.4$ ($\textsc{nli} = 16$) with $\eps = 10^{-2}$.\label{tab::resultscd1}}
\end{table}

We now provide a numerical insight on the comparison between the proposed IPM approach
and the commonly used semismooth Newton approach \cite{HIK02}.
We therefore compare \ipmbt and the implementation \ssnip of the global semismooth Newton method proposed for PDE-constrained optimization problems with sparsity-promoting terms in \cite{pss17}. When using the \ssnip approach, global convergence is attained using a nonsmooth line-search strategy
and the linear systems arising in the linear algebra phase are solved by 
using preconditioned {\scshape gmres}. We consider the  $2\times2$ block formulation and 
an indefinite preconditioner available in a factorized form \cite{pss17,pst15}. 
Since the semismooth approach requires a diagonal mass matrix in the discretization of the complementarity
conditions, in the experiments with \ssnip we use a lumped mass matrix.
Table \ref{tab::resultscd1_new} collects  results concerning the nonlinear behaviour of
the two methods: the number of nonlinear iterations ({\sc nli}) and the total CPU time ({\sc tcpu}).

It is interesting to note that the number of nonlinear Interior Point iterations does not vary with $\alpha$.
In fact, the mildly aggressive choice of barrier reduction factor $\sigma$ yields  a low number of nonlinear iterations,
 even for limiting values of $\alpha$.
By contrast, \ssnip  struggles as $\alpha \rightarrow 0$. Furthermore, overall the 
Interior Point strategy outperforms the semismooth method in terms of total CPU time.\\


\begin{table}[htb!]
\begin{center}
\begin{tabular}{llrrrr}
\toprule
 &                               &      \multicolumn{ 2}{c}{\ipmbt}     &     \multicolumn{ 2}{|c }{\ssnip}\\
\midrule
  $h=2^{-\ell}$ &   $\mathrm{log}_{10}\alpha$ &  {\sc nli}  &  {\sc tcpu} & {\sc nli}  &  {\sc tcpu}  \\
\midrule
       
\multicolumn{ 1}{c}{6} &         -2 &         11 &        2.8 &          5 &        4.2  \\

\multicolumn{ 1}{c}{} &         -4 &         11 &        2.5 &         19 &       27.9 \\

\multicolumn{ 1}{c}{} &         -6 &         11 &        2.4 &       $>100$ &                       \\

\multicolumn{ 1}{c}{} &         -8 &         11 &        2.4 &       $>100$ &                     \\

\multicolumn{ 1}{c}{7} &         -2 &         11 &        9.4 &          5 &       14.0     \\

\multicolumn{ 1}{c}{} &         -4 &         11 &        8.7 &         18 &      101.9     \\

\multicolumn{ 1}{c}{} &         -6 &         11 &        8.7 &       $>100$ &                      \\

\multicolumn{ 1}{c}{} &         -8 &         11 &        9.1 &       $>100$ &               \\

\multicolumn{ 1}{c}{8} &         -2 &         11 &       36.6 &          5 &       43.4      \\

\multicolumn{ 1}{c}{} &         -4 &         11 &       34.4 &         20 &      345.3    \\

\multicolumn{ 1}{c}{} &         -6 &         11 &       33.9 &       $>100$ &                   \\

\multicolumn{ 1}{c}{} &         -8 &         11 &       33.8 &       $>100$ &                     \\

\multicolumn{ 1}{c}{9} &         -2 &         11 &      155.9 &          5 &      147.3           \\

\multicolumn{ 1}{c}{} &         -4 &         11 &      149.8 &         21 &     1265.4            \\

\multicolumn{ 1}{c}{} &         -6 &         11 &      148.9 &       $>100$ &           \\

\multicolumn{ 1}{c}{} &         -8 &         11 &      149.6 &       $>100$ &                \\

\bottomrule
\end{tabular}  
\end{center}
\caption{Convection--diffusion problem: comparison between \ipmbt and \ssnip in terms of nonlinear iterations and total CPU times for problem with control constraints, for a range of $h$ and $\alpha$,
 $\beta = 10^{-3}$, $\epsilon = 10^{-1}$.  \label{tab::resultscd1_new}}
\end{table}

\subsection{A Heat Equation Problem}
To demonstrate the applicability of our methodology to time-dependent problems, we now perform experiments on an optimization problem with the heat equation acting as a constraint. We utilize the implicit Euler scheme on a time interval up to $T=1$, for varying values of time-step $\tau$, and set a time-independent desired state to be $\rm y_d=\sin(\pi {\rm x_1})\sin(\pi {\rm x_2}) $. We consider a control problem with full observations, with Table \ref{tab::resultsheat1} illustrating the performance of the Interior Point method and preconditioner $\mathcal{P}_T$ for varying mesh-sizes and values of $\alpha$, with fixed $\beta=10^{-2}$. Considerable robustness is again achieved, in particular with respect to changes in the time-step.
\begin{table}[htb!]
\begin{center}
\begin{tabular}{llrrrrrr}
\toprule
           &            & \multicolumn{ 6}{c}{\ipmbt} \\
					\midrule
           &            & \multicolumn{ 2}{c}{$\tau = 0.04$ } & \multicolumn{ 2}{c}{$\tau = 0.02$ } & \multicolumn{ 2}{c}{$\tau = 0.01$ } \\
					\midrule
       $h=2^{-\ell}$ &   $\mathrm{log}_{10}\alpha$ &  {\sc av-li}  &  {\sc av-cpu} &      {\sc av-li}  &  {\sc av-cpu} &  {\sc av-li}  &  {\sc av-cpu} \\
\midrule
\multicolumn{ 1}{c}{4} &         $-2$ &       13.9 &        0.6 &       13.1 &        1.0 &       13.1 &        2.2 \\

\multicolumn{ 1}{c}{} &         $-4$ &       13.3 &        0.5 &       12.2 &        1.0 &       12.3 &        2.0 \\

\multicolumn{ 1}{c}{} &         $-6$ &       12.8 &        0.5 &       12.0 &        1.0 &       12.0 &        2.0 \\

\multicolumn{ 1}{c}{5} &         $-2$ &       14.6 &        1.6 &       14.0 &        3.1 &       14.7 &        6.6 \\

\multicolumn{ 1}{c}{} &         $-4$ &       13.9 &        1.5 &       13.3 &        2.9 &       13.3 &        5.8 \\

\multicolumn{ 1}{c}{} &         $-6$ &       13.6 &        1.5 &       12.8 &        2.8 &       13.0 &        5.7 \\

\multicolumn{ 1}{c}{6} &         $-2$ &       15.5 &        5.9 &       14.6 &       11.4 &       15.4 &       23.7 \\

\multicolumn{ 1}{c}{} &         $-4$ &       14.8 &        5.8 &       14.0 &       10.6 &       14.0 &       21.7 \\

\multicolumn{ 1}{c}{} &         $-6$ &       14.6 &        5.5 &       13.8 &       10.6 &       13.9 &       21.5 \\
\bottomrule
\end{tabular}  
\end{center}
\caption{Heat equation problem: average Krylov iterations and CPU times for problem with control constraints, 
for a range of $h$, $\alpha$, and $\tau$, $\beta = 10^{-2}$, $\sigma=0.25$ ($\textsc{nli} = 13$). 
\label{tab::resultsheat1}}
\end{table}

\vspace{1em}

\begin{remark} We highlight that the number of nonlinear Interior Point iterations almost does not vary with $\alpha$, due
to the suitable choices made for the barrier reduction factor $\sigma$. In particular, in all the test cases
discussed, the choice of $\sigma$ is mildly aggressive (from $0.2$ to $0.4$ in the most difficult cases),
yielding a low number of nonlinear iterations, even for limiting values of $\alpha$.
By contrast,  a semismooth Newton approach globalized with a line-search
strategy may perform poorly as $\alpha \rightarrow 0$.
\end{remark}

\section{Conclusions}

We have presented a new Interior Point method for PDE-constrained optimization problems that include additional box constraints on the control variable, as well as possibly the state variable, and a sparsity-promoting $\rm L^1$-norm term for the control within the cost functional. We incorporated a splitting of the control into positive and negative parts, as well as a suitable nodal quadrature rule, to linearize the $\rm L^1$-norm, and considered preconditioned iterative solvers for the Newton systems arising at each Interior Point iteration. Through theoretical justification for our approximations of the $(1,1)$-block and Schur complement of the Newton systems, as well as numerical experiments, we have demonstrated the effectiveness and robustness of our approach, which may be applied within symmetric and non-symmetric Krylov methods, for a range of steady and time-dependent PDE-constrained optimization problems.

\Appendix
\section{Interior Point Algorithm for Quadratic Programming}\label{IPalgo}
In the Algorithm below, we present the structure of the Interior Point method that we apply within our numerical experiments, following the Interior Point path-following scheme described in \cite{gondzio12}. It is clear that the main computational effort arises from solving the Newton system \eqref{NewtonSystem} at each iteration.

\algo{ipm_algo}{Interior Point Algorithm for Quadratic Programming}{\vspace{-2em}
\begin{align*}
\ &\textbf{Parameters} \\
\ &\quad\quad\alpha_{0} \in(0,1),~~\text{step-size factor to boundary} \\
\ &\quad\quad\sigma\in(0,1),~~\text{barrier reduction parameter} \\
\ &\quad\quad\epsilon_{p},~\epsilon_{d},~\epsilon_{c},~~\text{stopping tolerances} \\
\ &\quad\quad\text{Interior point method stops when }\big\|{\xi}_{p}^{k}\big\|\leq\epsilon_{p},~\big\|{\xi}_{d}^{k}\big\|\leq\epsilon_{d},~\big\|{\xi}_{c}^{k}\big\|\leq\epsilon_{c} \\
\ &\textbf{Initialize IPM} \\
\ &\quad\quad\text{Set the initial guesses for }{y}^{0},~{z}^{0},~{p}^{0},~{\lambda}_{y,a}^{0},~{\lambda}_{y,b}^{0},~{\lambda}_{z,a}^{0},~{\lambda}_{z,b}^{0} \\
\ &\quad\quad\text{Set the initial barrier parameter }\mu^{0} \\
\ &\quad\quad\text{Compute primal infeasibility } {\xi}_{p}^{0}, \text{ dual infeasibility } {\xi}_{d}^{0}, \text{ and} 
\text{ complementarity gap }{\xi}_{c}^{0}, \\
\ &\quad\quad\quad\quad \text{as in \eqref{prdu}--\eqref{gap} with }k=0 \\ 
\ &\textbf{Interior Point Method} \\
\ &\quad\quad\text{while}~~\left(\big\|{\xi}_{p}^{k}\big\|>\epsilon_{p}~~\text{or}~~\big\|{\xi}_{d}^{k}\big\|>\epsilon_{d}~~\text{or}~~\big\|{\xi}_{c}^{k}\big\|>\epsilon_{c}\right) \\
\ &\quad\quad\quad\quad\text{Reduce barrier parameter}~\mu^{k+1}=\sigma\mu^{k} \\
\ &\quad\quad\quad\quad\text{Solve Newton system }\eqref{NewtonSystem}\text{ for primal-dual Newton direction}~{\Delta}{y},~{\Delta}{z},~{\Delta p} \\
\ &\quad\quad\quad\quad\text{Use }\text{\eqref{zupdate1}--\eqref{zupdate4}}\text{ to find }{\Delta}{\lambda}_{y,a},~{\Delta}{\lambda}_{y,b},~{\Delta}{\lambda}_{z,a},~{\Delta}{\lambda}_{z,b} \\
\ &\quad\quad\quad\quad\text{Find }\alpha_{P},~\alpha_{D}~\text{s.t. bound constraints on primal and dual variables hold} \\
\ &\quad\quad\quad\quad\text{Set }\alpha_{P}=\alpha_{0}\alpha_{P},~\alpha_{D}=\alpha_{0}\alpha_{D} \\
\ &\quad\quad\quad\quad\text{Make step: }{y}^{k+1}={y}^{k}+\alpha_{P}{\Delta}{y},~{z}^{k+1}={z}^{k}+\alpha_{P}{\Delta}{z},~{p}^{k+1}={p}^{k}+\alpha_{D}{\Delta p} \\
\ &\quad\quad\quad\quad\quad\quad{\lambda}_{y,a}^{k+1}={\lambda}_{y,a}^{k}+\alpha_{D}{\Delta}{\lambda}_{y,a}, \ 
                                 {\lambda}_{y,b}^{k+1}={\lambda}_{y,b}^{k}+\alpha_{D}{\Delta}{\lambda}_{y,b} \\
\ &\quad\quad\quad\quad\quad\quad{\lambda}_{z,a}^{k+1}={\lambda}_{z,a}^{k}+\alpha_{D}{\Delta}{\lambda}_{z,a}, \ 
                                 {\lambda}_{z,b}^{k+1}={\lambda}_{z,b}^{k}+\alpha_{D}{\Delta}{\lambda}_{z,b} \\
\ &\quad\quad\quad\quad\text{Update infeasibilities }  {\xi}_{p}^{k+1},~{\xi}_{d}^{k+1}, \text{ and compute the complementarity gap }  {\xi}_{c}^{k+1} \\
\ &\quad\quad\quad\quad\quad\quad\text{as in \eqref{prdu}--\eqref{gap}} \\
%
%
\ &\quad\quad\quad\quad\text{Set iteration number }k=k+1 \\
\ &\quad\quad\text{end}
\end{align*}\vspace{-1.5em}
}

\textbf{Acknowledgments.} 
J. W. Pearson gratefully acknowledges support from the Engineering and Physical Sciences Research Council (EPSRC) Fellowship EP/M018857/2, and a Fellowship from The Alan Turing Institute in London.
M. Porcelli and M. Stoll were partially supported by the {\em DAAD-MIUR Joint Mobility Program} 2018--2020 (Grant 57396654).
The work of M. Porcelli was also partially supported by the {\em National Group of Computing Science (GNCS-INDAM)}.

\bibliographystyle{siam}
\bibliography{data}

\begin{thebibliography}{10}

\bibitem{dealii}
{\sc W.~Bangerth, R.~Hartmann, and G.~Kanschat}, {\em deal.{II}---{A}
  general-purpose object-oriented finite element library}, ACM Transactions on
  Mathematical Software, 33 (2007), p.~Art. 24.

\bibitem{BenDOS15}
{\sc P.~Benner, S.~Dolgov, A.~Onwunta, and M.~Stoll}, {\em {Low-rank solvers
  for unsteady Stokes--Brinkman optimal control problem with random data}},
  Computer Methods in Applied Mechanics and Engineering, 304 (2016),
  pp.~26--54.

\bibitem{BenGolLie05}
{\sc M.~Benzi, G.~H. Golub, and J.~Liesen}, {\em Numerical solution of saddle
  point problems}, Acta Numerica, 14 (2005), pp.~1--137.

\bibitem{BIK99}
{\sc M.~Bergounioux, K.~Ito, and K.~Kunisch}, {\em Primal-dual strategy for
  constrained optimal control problems}, SIAM Journal on Control and
  Optimization, 37 (1999), pp.~1176--1194.

\bibitem{Boyle2007}
{\sc J.~Boyle, M.~D. Mihajlovi\'c, and J.~A. Scott}, {\em {HSL\_MI20: An
  efficient AMG preconditioner for finite element problems in 3D}},
  International Journal for Numerical Methods in Engineering, 82 (2010),
  pp.~64--98.

\bibitem{BroH82}
{\sc A.~N. Brooks and T.~J.~R. Hughes}, {\em Streamline
  upwind/{Petrov--Galerkin} formulations for convection dominated flows with
  particular emphasis on the incompressible {Navier--Stokes} equations},
  Computer Methods in Applied Mechanics and Engineering, 32 (1982),
  pp.~199--259.

\bibitem{Cas86}
{\sc E.~Casas}, {\em Control of an elliptic problem with pointwise state
  constraints}, SIAM Journal on Control and Optimization, 24 (1986),
  pp.~1309--1318.

\bibitem{de2013image}
{\sc J.~C. De~los Reyes and C.-B. Sch{\"o}nlieb}, {\em Image denoising:
  Learning the noise model via nonsmooth {PDE}-constrained optimization},
  Inverse Problems \& Imaging, 7 (2013), pp.~1183--1214.

\bibitem{ifissmatlab}
{\sc H.~C. Elman, A.~Ramage, and D.~J. Silvester}, {\em Algorithm 866: {IFISS},
  a {M}atlab toolbox for modelling incompressible flow}, ACM Transactions on
  Mathematical Software, 33 (2007), p.~Art. 14.

\bibitem{ifisslink}
\leavevmode\vrule height 2pt depth -1.6pt width 23pt, {\em Incompressible
  {F}low and {I}terative {S}olver {S}oftware ({IFISS})}, Version 3.5,
  \texttt{http://www.maths.manchester.ac.uk/}$\sim$\texttt{djs/ifiss/},
  (2018).

\bibitem{GPSR}
{\sc M.~A.~T. Figueiredo, R.~D. Nowak, and S.~J. Wright}, {\em Gradient
  projection for sparse reconstruction: Application to compressed sensing and
  other inverse problems}, IEEE Journal of Selected Topics in Signal
  Processing, 1 (2007), pp.~586--597.

\bibitem{FG-pseudo16}
{\sc K.~Fountoulakis and J.~Gondzio}, {\em A second-order method for strongly
  convex $\ell_1$-regularization problems}, Mathematical Programming, 156
  (2016), pp.~189--219.

\bibitem{FG-IPM14}
{\sc K.~Fountoulakis, J.~Gondzio, and P.~Zhlobich}, {\em Matrix-free interior
  point method for compressed sensing problems}, Mathematical Programming
  Computation, 6 (2014), pp.~1--31.

\bibitem{VGI61}
{\sc G.~H. Golub and R.~S. Varga}, {\em Chebyshev semi-iterative methods,
  successive over-relaxation iterative methods, and second order {R}ichardson
  iterative methods. {I}}, Numerische Mathematik, 3 (1961), pp.~147--156.

\bibitem{VGII61}
\leavevmode\vrule height 2pt depth -1.6pt width 23pt, {\em Chebyshev
  semi-iterative methods, successive over-relaxation iterative methods, and
  second order {R}ichardson iterative methods. {II}}, Numerische Mathematik, 3
  (1961), pp.~157--168.

\bibitem{gondzio12}
{\sc J.~Gondzio}, {\em Interior point methods 25 years later}, European Journal
  of Operational Research, 218 (2012), pp.~587--601.

\bibitem{gunther2012posteriori}
{\sc A.~G{\"u}nther, M.~Hinze, and M.~H. Tber}, {\em A posteriori error
  representations for elliptic optimal control problems with control and state
  constraints}, in Constrained Optimization and Optimal Control for Partial
  Differential Equations, Springer, 2012, pp.~303--317.

\bibitem{HerOW15}
{\sc R.~Herzog, J.~Obermeier, and G.~Wachsmuth}, {\em Annular and sectorial
  sparsity in optimal control of elliptic equations}, Computational
  Optimization and Applications, 62 (2015), pp.~157--180.

\bibitem{herzog2018fast}
{\sc R.~Herzog, J.~W. Pearson, and M.~Stoll}, {\em Fast iterative solvers for
  an optimal transport problem}, Advances in Computational Mathematics,
  doi:10.1007/s10444-018-9625-5 (2018).

\bibitem{HS10}
{\sc R.~Herzog and E.~W. Sachs}, {\em Preconditioned conjugate gradient method
  for optimal control problems with control and state constraints}, SIAM
  Journal on Matrix Analysis and Applications, 31 (2010), pp.~2291--2317.

\bibitem{HSW11_DS}
{\sc R.~Herzog, G.~Stadler, and G.~Wachsmuth}, {\em {Directional sparsity in
  optimal control of partial differential equations}}, {SIAM Journal on Control
  and Optimization}, 50 (2012), pp.~943--963.

\bibitem{HIK02}
{\sc M.~Hinterm{\"u}ller, K.~Ito, and K.~Kunisch}, {\em The primal-dual active
  set strategy as a semismooth {N}ewton method}, SIAM Journal on Optimization,
  13 (2002), pp.~865--888.

\bibitem{Hinze2000}
{\sc M.~Hinze}, {\em Optimal and instantaneous control of the instationary
  Navier--Stokes equations}, Habilitation, Technische Universit\"{a}t Berlin,
  2000.

\bibitem{book::hpuu09}
{\sc M.~Hinze, R.~Pinnau, M.~Ulbrich, and S.~Ulbrich}, {\em Optimization with
  {PDE} Constraints}, Mathematical Modelling: Theory and Applications,
  Springer-Verlag, New York, 2009.

\bibitem{Ipsen}
{\sc I.~Ipsen}, {\em A note on preconditioning non-symmetric matrices}, SIAM
  Journal on Scientific Computing, 23 (2001), pp.~1050--1051.

\bibitem{book::IK08}
{\sc K.~Ito and K.~Kunisch}, {\em Lagrange Multiplier Approach to Variational
  Problems and Applications}, vol.~15 of Advances in Design and Control,
  Society for Industrial and Applied Mathematics, Philadelphia, PA, 2008.

\bibitem{Ku95}
{\sc Y.~A. Kuznetsov}, {\em Efficient iterative solvers for elliptic finite
  element problems on nonmatching grids}, Russian Journal of Numerical Analysis
  and Mathematical Modelling, 10 (1995), pp.~187--211.

\bibitem{LSinverses02}
{\sc T.-T. Lu and S.-H. Shiou}, {\em Inverses of $2\times2$ block matrices},
  Computers \& Mathematics with Applications, 43 (2002), pp.~119--129.

\bibitem{preconMGW}
{\sc M.~F. Murphy, G.~H. Golub, and A.~J. Wathen}, {\em A note on
  preconditioning for indefinite linear systems}, SIAM Journal on Scientific
  Computing, 21 (2000), pp.~1969--1972.

\bibitem{NocW06}
{\sc J.~Nocedal and S.~J. Wright}, {\em Numerical Optimization}, Springer
  Series in Operations Research and Financial Engineering, Springer, New York,
  2nd~ed., 2006.

\bibitem{minres}
{\sc C.~C. Paige and M.~A. Saunders}, {\em Solution of sparse indefinite
  systems of linear equations}, SIAM Journal on Numerical Analysis, 12 (1975),
  pp.~617--629.

\bibitem{PGIP17}
{\sc J.~W. Pearson and J.~Gondzio}, {\em Fast interior point solution of
  quadratic programming problems arising from {PDE}-constrained optimization},
  Numerische Mathematik, 137 (2017), pp.~959--999.

\bibitem{PSW11}
{\sc J.~W. Pearson, M.~Stoll, and A.~J. Wathen}, {\em {Regularization-robust
  preconditioners for time-dependent {PDE}-constrained optimization problems}},
  SIAM Journal on Matrix Analysis and Applications, 33 (2012), pp.~1126--1152.

\bibitem{PW10}
{\sc J.~W. Pearson and A.~J. Wathen}, {\em {A new approximation of the Schur
  complement in preconditioners for PDE-constrained optimization}}, Numerical
  Linear Algebra with Applications, 19 (2012), pp.~816--829.

\bibitem{PW11}
\leavevmode\vrule height 2pt depth -1.6pt width 23pt, {\em Fast iterative
  solvers for convection--diffusion control problems}, Electronic Transactions
  on Numerical Analysis, 40 (2013), pp.~294--310.

\bibitem{pss17}
{\sc M.~Porcelli, V.~Simoncini, and M.~Stoll}, {\em Preconditioning
  {PDE}-constrained optimization with ${L}^1$-sparsity and control
  constraints}, Computers \& Mathematics with Applications, 74 (2017),
  pp.~1059--1075.

\bibitem{pst15}
{\sc M.~Porcelli, V.~Simoncini, and M.~Tani}, {\em Preconditioning of
  active-set {N}ewton methods for {PDE}-constrained optimal control problems},
  SIAM Journal on Scientific Computing, 37 (2015), pp.~S472--S502.

\bibitem{RSW09}
{\sc T.~Rees, M.~Stoll, and A.~Wathen}, {\em {All-at-once preconditioners for
  PDE-constrained optimization}}, Kybernetika, 46 (2010), pp.~341--360.

\bibitem{gmres}
{\sc Y.~Saad and M.~H. Schultz}, {\em G{MRES}: {A} generalized minimal residual
  algorithm for solving nonsymmetric linear systems}, SIAM Journal on
  Scientific and Statistical Computing, 7 (1986), pp.~856--869.

\bibitem{SongYu3}
{\sc X.~Song, B.~Chen, and B.~Yu}, {\em Error estimates for sparse optimal
  control problems by piecewise linear finite element approximation}, arXiv
  preprint arXiv:1709.09539,  (2017).

\bibitem{SongYu2}
\leavevmode\vrule height 2pt depth -1.6pt width 23pt, {\em Mesh independence of
  an accelerated block coordinate descent method for sparse optimal control
  problems}, arXiv preprint arXiv:1709.00005,  (2017).

\bibitem{SongYu1}
\leavevmode\vrule height 2pt depth -1.6pt width 23pt, {\em An efficient
  duality-based approach for {PDE}-constrained sparse optimization},
  Computational Optimization and Applications, 69 (2018), pp.~461--500.

\bibitem{Sta09}
{\sc G.~Stadler}, {\em {Elliptic optimal control problems with $L^1$-control
  cost and applications for the placement of control devices}}, {Computational
  Optimization and Applications}, 44 (2009), pp.~159--181.

\bibitem{book::FT2010}
{\sc F.~Tr{\"o}ltzsch}, {\em {Optimal Control of Partial Differential
  Equations: Theory, Methods and Applications}}, American Mathematical Society,
  2010.

\bibitem{ulbrich2009primal}
{\sc M.~Ulbrich and S.~Ulbrich}, {\em Primal-dual interior-point methods for
  {PDE}-constrained optimization}, Mathematical Programming, 117 (2009),
  pp.~435--485.

\bibitem{wwL1}
{\sc G.~Wachsmuth and D.~Wachsmuth}, {\em Convergence and regularization
  results for optimal control problems with sparsity functional}, ESAIM:
  Control, Optimisation and Calculus of Variations, 17 (2011), pp.~858--886.

\bibitem{wathen87}
{\sc A.~J. Wathen}, {\em Realistic eigenvalue bounds for the {G}alerkin mass
  matrix}, IMA Journal of Numerical Analysis, 7 (1987), pp.~449--457.

\bibitem{RW08}
{\sc A.~J. Wathen and T.~Rees}, {\em Chebyshev semi-iteration in
  preconditioning for problems including the mass matrix}, Electronic
  Transactions in Numerical Analysis, 34 (2008), pp.~125--135.

\bibitem{fpcas}
{\sc Z.~Wen, W.~Yin, D.~Goldfarb, and Y.~Zhang}, {\em A fast algorithm for
  sparse reconstruction based on shrinkage, subspace optimization, and
  continuation}, SIAM Journal on Scientific Computing, 32 (2010),
  pp.~1832--1857.

\bibitem{IPMWright}
{\sc S.~J. Wright}, {\em Primal-Dual Interior-Point Methods}, Society for
  Industrial and Applied Mathematics, Philadelphia, PA, 1997.

\end{thebibliography}
\end{document}